\documentclass[titlepage]{article}
\usepackage{graphicx} 
\usepackage[margin=1in]{geometry}
\usepackage[utf8]{inputenc}
\usepackage{amsmath}
\usepackage{amsthm}
\usepackage{amsfonts}
\usepackage{amssymb}
\usepackage{ragged2e}
\usepackage{authblk}
\usepackage[citestyle=numeric-comp,]{biblatex}
\usepackage{hyperref}

\providecommand{\keywords}[1]
{
  \small	
  \textbf{\textit{Keywords: }} #1
}

\providecommand{\subjclass}[1]
{
  \small	
  \textbf{\textit{Mathematics Subject Classification: }} #1
}

\newcommand{\gaussian}{\genfrac{[}{]}{0pt}{}}
\newtheorem{theorem}{Theorem}
\newtheorem{corollary}[theorem]{Corollary}
\newtheorem{lemma}[theorem]{Lemma}
\newtheorem{proposition}[theorem]{Proposition}
\title{On certain $q$-multiple sums}
\author{\emph{Aung Phone Maw}} 
\date{September 2024}
\begin{document}
\maketitle

\raggedright\textbf{\Large{Abstract}}\\ \vspace{0.1cm}
    \quad We present outlines of a general method to reach certain kinds of $q$-multiple sum identities. Throughout our exposition, we shall give generalizations to the results given by Dilcher, Prodinger, Fu and Lascoux, Zeng, and Guo and Zhang concerning $q$-series identities related to divisor functions. Our exposition shall also provide a generalization of the duality relation for finite multiple harmonic $q$-series given by Bradley. Utilizing these generalizations, we will also arrive at some new interesting classes of $q$-multiple sums.

    \vspace{0.1cm}
    
\raggedright 
. \\
.
\par
\section{Introduction}

Throughout this paper, we shall use the following standard notation

\[ (x;q)_\infty = (x)_{\infty} = \prod\limits_{i \geq 1} (1-xq^{i-1}),
\]
\[
(x;q)_n = \prod\limits_{1 \leq i \leq n} (1-xq^{i-1}), \hspace{0.05cm} n \in \mathbb{N},
\]
\[ \hspace{0.05cm}  (x;q)_0 = 1.
\]

The study of \emph{$q$-identities related to divisor functions} \cite{KDilcher,FuLascoux,GuoZhang,GuoZeng,ZuCen,Xu,prodinger1999qanalogueformulahernandezobtained,prodinger2004qidentitiesfulascouxproved,Zeng}, has given rise to numerous interesting $q$-\\combinatorial identities. In these studies, the regular appearance of multiple sums is of noticeable significance. We shall now point out those $q$-combinatorial identities with multiple sums which will be examined throughout our study. The first appearance of these identities occurs in \cite{KDilcher}, where Dilcher gave the following identity, which holds for $m \geq 1$.

\[ \tag{1.1} \label{dilcher}
\sum\limits_{n \geq i_1 \geq ... \geq i_m \geq 1} \frac{q^{i_1+...+i_m}}{(1-q^{i_1})...(1-q^{i_m})} = \sum\limits_{1 \leq r \leq n} \gaussian{n}{r}\frac{(-1)^{r-1}q^{\binom{r}{2}+rm}}{(1-q^r)^m},
\]

where 

\[ \gaussian{n}{r} = \frac{[n][n-1]...[n-r+1]}{[r][r-1]...[1]}, \]

is the \emph{Gaussian binomial coefficient}, and $[n] = \frac{1-q^n}{1-q}$ is the \emph{q-number}. \eqref{dilcher} is provided as a certain analogue of the identity 

\[ \tag{1.2}
\sum\limits_{i \geq 1} \frac{q^i}{1-q^i} = \sum\limits_{r \geq 1} \frac{(-1)^{r-1}q^{\binom{r+1}{2}}}{(q;q)_r (1-q^r)}.
\]

given in \cite{Kluyver}. Note that the series on the left-hand side generates the arithmetic function of the number of divisors of a given natural number. Later, Prodinger \cite{prodinger1999qanalogueformulahernandezobtained} proved the following 

\[ \tag{1.3} \label{prodinger}
\sum\limits_{1 \leq i \leq n} \frac{q^{i(m-1)}}{(1-q^i)^m} = \sum\limits_{1 \leq r \leq n} \gaussian{n}{r}(-1)^{r-1}q^{\binom{r}{2}-rn} \sum\limits_{r = i_1 \geq ... \geq i_m \geq 1} \frac{q^{i_1+...+i_m}}{(1-q^{i_1})...(1-q^{i_m})},
\]

by inverting the original result of Dilcher and thus giving a $q$-analog of a formula of Hernández \cite{VHernandez}. Later, Fu and Lascoux \cite{FuLascoux} further generalized \eqref{dilcher} as 
\[ \tag{1.4} \label{FuLascoux}
\sum\limits_{n \geq i_1 \geq ... \geq i_m \geq 1} (-1)^{i_m-1}(x^{i_m}-(-1)^{i_m})\frac{q^{i_1+...+i_m}}{(1-q^{i_1})...(1-q^{i_m})} = \sum\limits_{1 \leq r \leq n} \gaussian{n}{r}\frac{(-1)^{r-1} x^{r} (-x^{-1};q)_rq^{rm}}{(1-q^r)^m}.
\]

For this, Fu and Lascoux used the Newton interpolation. Meanwhile, Prodinger \cite{prodinger2004qidentitiesfulascouxproved} and Zeng \cite{Zeng} provided different proofs of \eqref{FuLascoux}. Zeng \cite{Zeng} in particular, using the method of partial fraction decomposition, obtained a further generalization of \eqref{FuLascoux}, which can be stated as 

\[ \tag{1.5} \label{Zeng}
\frac{(q;q)_n}{(z;q)_n}\sum\limits_{n \geq i_1 \geq ... \geq i_m \geq 1} \frac{x^{i_m}(zq;q)_{i_m}}{(q;q)_{i_m}}\frac{q^{i_1+...+i_m}}{(1-zq^{i_1})...(1-zq^{i_m})} = \sum\limits_{1 \leq r \leq n} \gaussian{n}{r}\frac{(x^r(x^{-1};q)_r+(-1)^{r-1}q^{\binom{r}{2}})q^{rm}}{(1-zq^r)^m}.
\]
This is a very common generalization of \eqref{dilcher}. Ismail and Stanton \cite{IshStanton}, and A. Xu \cite{Xu} also provided different proofs of \eqref{Zeng}. Furthermore, we note that Guo and Zhang \cite{GuoZhang} obtained the following very unique generalization of \eqref{dilcher} 
\[ \tag{1.6} \label{gouzhang}
-\sum\limits_{n \geq i_1 \geq ... \geq i_m \geq 1} \frac{q^{i_1+i_2...+i_m}}{(1-q^{i_1})(1-q^{i_2})...(1-q^{i_m})(1-zq^{i_1-1})(1-zq^{i_2-2})...(1-zq^{i_m-m})} = \sum\limits_{1 \leq r \leq n} \gaussian{n}{r} \frac{(z^{-1}q^m;q)_r(z;q)_{n-r}}{(zq^{-m};q)_{m+n}(1-q^r)^m}z^r.
\]

Lastly we also add that, in fact, \eqref{dilcher} and \eqref{prodinger} can be viewed as particular cases of the duality identity (Theorem 1) appearing in Bradley \cite{Bradley}. Before stating our main results, let us define the following \emph{generalized q-multiple harmonic sums}.

\[
H_n[s_1,s_2,...,s_k;x:q] = H_n[s_1,s_2,...,s_k;x] = \frac{q^{s_1-1}}{[n]^{s_1-1}} \sum\limits_{n \geq i_1 \geq ... \geq i_{k-1} \geq 1} \frac{q^{i_1(s_2-1)+...+i_{k-1}(s_k-1)}}{[i_1]^{s_2}...[i_{k-1}]^{s_k}} x^{i_{k-1}},
\]

\[
U_n[s_1,s_2,...,s_k;x:q] = U_n[s_1,s_2,...,s_k;x] = \frac{1}{[n]^{s_1-1}} \sum\limits_{n \geq i_1 \geq ... \geq i_{k-1} \geq 1} 
\frac{q^{i_1+...+i_{k-1}}}{[i_1]^{s_2}...[i_{k-1}]^{s_k}} (1-(x;q)_{i_{k-1}}).
\]

For $k$ copies of the argument $n$ we shall write $\{n\}^k$. For example, $H_n[\{1\}^3, 2;x] = H_n[1,1,1,2;x]$, $H_n[3, \{2\}^4, 5;x]= H_n[3,2,2,2,2,5;x]$, and $U_n[2,\{1\}^4;x] = U_n[2,1,1,1,1;x]$. Then our first result can be stated as follows.
\vspace{0.7cm}
\begin{theorem} Let $A_1(q), A_2(q), ...$ and $B_1(q), B_2(q), ...$ be sequences not depending on $x$, satisfying the relation 

\[
\sum\limits_{r\geq 1} A_r(q) x^r = \sum\limits_{r \geq 1} B_r(q) (1-(x;q)_r),
\]

for all complex values $x$, then

\begin{multline} \notag
\sum\limits_{r \geq 1} A_r(q) H_r[m_1+1, \{1\}^{n_1-1}, m_2+1, \{1\}^{n_2-1}, m_3+1,..., \{1\}^{n_{k-1}-1}, m_k+1, \{1\}^{n_k} ; x] \\= \sum\limits_{r \geq 1} B_r(q) U_r[\{1\}^{m_1}, n_1+1, \{1\}^{m_2-1}, n_2+1, \{1\}^{m_3-1}, n_3+1,..., \{1\}^{m_k-1},n_k+1;x],
\end{multline}
for all non-negative integers, $m_1, m_2,...,m_k$ and $n_1,n_2,...,n_k$. Where it is understood that when $s=0$, $H_N[n_1,...,n_p, n+1, \{1\}^{s-1}, r+1,r_1,...,r_q;x] = H_N[n_1,...,n_p, n+r+1,r_1,...,r_q;x]$ and $U_N[n_1,...,n_p, n+1, \{1\}^{s-1}, r+1,r_1,...,r_q;x] = U_N[n_1,...,n_p, n+r+1,r_1,...,r_q;x]$. 
\end{theorem}
It will be evident that \textbf{Theorem 1} encompasses Bradley's duality relation and the identities \eqref{dilcher}-\eqref{FuLascoux}. Next, we shall state a further analog.
\vspace{0.7cm}

\begin{theorem} Let $A_1(q), A_2(q), ...$ and $B_1(q), B_2(q), ...$ be sequences not depending on $x$, satisfying the relation 

\[
\sum\limits_{r\geq 1} A_r(q) x^r = \sum\limits_{r \geq 1} B_r(q) (1-(x;q)_r),
\]

for all complex values $x$, then

\begin{multline} \notag
    \sum\limits_{r\geq 1} A_r(q)\frac{q^r(q;q)_r}{(1-z_1q^r)(y_1q;q)_r}\sum\limits_{r_0=r \geq r_1 \geq r_2 \geq \cdots \geq r_{k-1} \geq r_k \geq 1}  \frac{x^{r_{k}}y_{k}^{r_{k-1}-r_{k}}(y_kq;q)_{r_k-1}}{(q;q)_{r_k}}  \prod\limits_{j=1}^{k-1} \frac{q^{r_j}y_j^{r_{j-1}-r_j}(y_jq;q)_{r_j-1}}{(1-z_{j+1}q^{r_j})(y_{j+1}q;q)_{r_j}} 
\end{multline}

\begin{multline} \notag
= \sum\limits_{r\geq 1} B_r(q) \frac{(q;q)_r}{(z_1q;q)_r}\sum\limits_{r_0=r \geq r_1 \geq r_2 \geq \cdots \geq r_{k-1} \geq r_{k} \geq 1} \frac{(1-(x;q)_{r_k})q^{r_{k}}(z_{k}q;q)_{r_{k}-1}}{(1-y_{k}q^{r_{k}})(q;q)_{r_{k}}}  \prod\limits_{j=1}^{k-1} \frac{q^{r_j}(z_jq;q)_{r_j-1}}{(1-y_jq^{r_j})(z_{j+1}q;q)_{r_j}} ,
\end{multline}

for all natural numbers $k$, and for all complex values $z_1,z_2,...,z_k$, and $y_1,y_2,...,y_k$, except at the points $q^{-r}, r \in \mathbb{N}$, where the expressions on both sides exhibit singularities.
    
\end{theorem}

It will be shown that \textbf{Theorem 2} encompasses the identities \eqref{Zeng} and \eqref{gouzhang}. The following proposition in particular is a generalization of Guo and Zhang's identities, \eqref{gouzhang} and Theorem 4.1 in \cite{GuoZhang}.
\vspace{0.7cm}
\begin{proposition} For natural numbers $n$ and $k$, and for all complex values $y,z,t$, there holds

    \begin{multline} \tag{1.7} \label{generalizedGouZhang}
 \sum\limits_{n \geq r_1 \geq r_2 \geq \cdots \geq r_{k-1} \geq r_k \geq 1} \left(\frac{(t;q)_{r_k}}{(ytq;q)_{r_k}}-\frac{(zy^{-1}q^{-k};q)_{r_k}}{(zq^{-k+1};q)_{r_k}}\right)\frac{q^{r_k}(yq;q)_{r_k-1}}{(q;q)_{r_k}} \prod\limits_{j=1}^{k-1}\frac{q^{r_j}}{(1-zq^{r_j-j})(1-yq^{r_j})}  \\= \frac{(yq;q)_n}{(q;q)_n(zq^{-k+1};q)_{k+n-1}}\sum\limits_{i \geq 1} \gaussian{n}{i} \frac{z^iy^{-i}(yq;q)_i(ytz^{-1}q^{k};q)_i(zy^{-1};q)_{n-i}}{(1-yq^i)^{k}(ytq;q)_i},
\end{multline}

except at the points $z \in \{ q^{k-1}, q^{k-2},..., q^{-n+1} \}$, $y \in \{q^{-1},...,q^{-n} \}$, and $t \in \{y^{-1}q^{-1},..., y^{-1}q^{-n}\}$, where the expressions from both sides exhibit singularities.
\end{proposition}
In deriving these statements, the following lemmas will prove to be useful for our manipulations of $q$-identities.
\vspace{0.7cm}
\begin{lemma}
    Let $A_1(q),A_2(q),...$ and $B_1(q), B_2(q),...$ be sequences not depending on $x$, satisfying the relation
    \[ \tag{1.8} \label{inverta}
    \sum\limits_{i \geq 1} A_i(q) x^i = \sum\limits_{i \geq 1} B_i(q) (1-(x;q)_i),
    \]
    for all complex values $x$, then 
    \[ \tag{1.9} \label{invertb}
    \sum\limits_{i \geq 1} A_i(q^{-1}) (1-(x;q)_i) = \sum\limits_{i \geq 1} B_i(q^{-1}) x^i.
    \]
\end{lemma}

\begin{proof}[Proof of Lemma 4]
Let us define the $q$-shift operator denoted as $\eta_x$ as defined in \cite{Andrews} 

\[
\eta_x f(x) = f(xq), \] 
\[ {\eta_x}^m f(x) = {\eta_x}^{m-1} \eta_x f(x), m \in \mathbb{N}, \] 
\[ {\eta_x}^0 f(x) = f(x).
\]

Then, let us consider applying the operator 
\[ \sum\limits_{m\geq 0} \frac{(x;q)_m}{(q;q)_m}y^m \eta_x^m, \] 
to both sides of \eqref{inverta}. Using the fact 
\[ (xq^m;q)_i= \frac{(xq^i;q)_m(x;q)_i}{(x;q)_m} ,\] 

and Heine's $q$-binomial theorem \cite{Heine} 

\[ \sum\limits_{m\geq0} \frac{(x;q)_m}{(q;q)_m}y^m = \frac{(xy;q)_{\infty}}{(y;q)_\infty}, \] we arrive at 

\[ \tag{1.10} \label{invertc}
    \sum\limits_{i \geq 1} A_i(q) x^i\frac{(y;q)_i}{(xy;q)_i} = \sum\limits_{i \geq 1} B_i(q) (1-\frac{(x;q)_i}{(xy;q)_i}).
    \]
Now we replace $q$ by $q^{-1}$, $x$ by $x^{-1}$, and $y$ by $y^{-1}$, to get

\[ \tag{1.11} \label{invertd}
    \sum\limits_{i \geq 1} A_i(q^{-1}) \frac{(y;q)_i}{(xy;q)_i} = \sum\limits_{i \geq 1} B_i(q^{-1}) (1-\frac{y^i(x;q)_i}{(xy;q)_i}).
    \]

Finally, we interchange $x$ and $y$, and put $y=0$ to arrive at \eqref{invertb}, and the proof is complete.
\end{proof}
We see that our proof of \textbf{Lemma 4} also implies the following lemma.

\vspace{0.7cm}
\begin{lemma}
    Let $A_1(q),A_2(q),...$ and $B_1(q), B_2(q),...$ be sequences not depending on $x$, satisfying the relation
    \[ 
    \sum\limits_{i \geq 1} A_i(q) x^i = \sum\limits_{i \geq 1} B_i(q) (1-(x;q)_i),
    \]
    for all complex values $x$, then 
    \[
    \sum\limits_{i \geq 1} A_i(q)x^i\frac{(y;q)_i}{(xy;q)_i}  = \sum\limits_{i \geq 1} B_i(q) (1-\frac{(x;q)_i}{(xy;q)_i}),
    \]
    for all complex values $x$ and $y$, except at the points $y = x^{-1}q^{-r}, r \in \mathbb{N} \cup \{0\}$, where the expressions from both sides exhibit singularities.
\end{lemma}

\textbf{Lemma 4} and \textbf{Lemma 5} will allow us to conveniently interchange between different forms of $q$-statements. In the last section, we shall provide a general transformation formula for certain types of basic hypergeometric multiple sums. In exploring some of its consequences we will also be able to provide a new class of $q-$multiple sums identities\\ such as 

\begin{multline}
   \frac{(q;q)_n}{(yt;q)_n-wt^n(y;q)_n}\left\{ 1+ w \sum\limits_{1 \leq r_1 < n} \frac{t^{r_1}(y;q)_{r_1}(t;q)_{n-r_1}}{((yt;q)_{r_1}-wt^{r_1}(y;q)_{r_1})(q;q)_{n-r_1}} \right. \\ \left. + w^2 \sum\limits_{1 \leq r_2 < r_1 < n} \frac{t^{r_1+r_2}(y;q)_{r_1}(t;q)_{n-r_1}(y;q)_{r_2}(t;q)_{r_1-r_2}}{((yt;q)_{r_1}-wt^{r_1}(y;q)_{r_1})(q;q)_{n-r_1}((yt;q)_{r_2}-wt^{r_2}(y;q)_{r_2})(q;q)_{r_1-r_2}} + ... \right. \\ \left. ...+ w^{n-1} \sum\limits_{1 \leq r_{n-1} <...<r_2<r_1<n} \ \prod\limits_{\substack{j=1 \\ r_0 = n}}^{n-1}\frac{t^{r_j}(y;q)_{r_j}(t;q)_{r_{j-1}-r_j}}{((yt;q)_{r_j}-wt^{r_j}(y;q)_{r_j})(q;q)_{r_{j-1}-r_j}}  \right\}  \\= \sum\limits_{r\geq 1} \gaussian{n}{r} \frac{(-1)^{r-1}q^{\binom{r}{2}}(1-q^r)(yt;q)_r}{(1-ytq^{r-1})((yt;q)_r-wt^r(y;q)_r)}, \tag{1.12} \label{egnewexpression1}
\end{multline}

and 

\begin{multline}
   \sum\limits_{r\geq 1} \gaussian{n}{r} (-1)^{r-1}q^{\binom{r+1}{2}-nr}\frac{(q;q)_r}{(zw;q)_r-tw^r(z;q)_r} \left\{1+ t \sum\limits_{1 \leq i_1 < r} \frac{w^{i_1}(z;q)_{i_1}(w;q)_{r-i_1}}{((zw;q)_{i_1}-tw^{i_1} (z;q)_{i_1})(q;q)_{r-i_1} } \right. \\ \left. + t^2\sum\limits_{1 \leq i_2 <i_1 < r} \frac{w^{i_1+i_2}(z;q)_{i_1}(w;q)_{r-i_1}(z;q)_{i_2}(w;q)_{i_1-i_2}}{((zw;q)_{i_1}-tw^{i_1} (z;q)_{i_1})(q;q)_{r-i_1}((zw;q)_{i_2}-tw^{i_2} (z;q)_{i_2})(q;q)_{i_1-i_2} } + ... \right. \\ \left. ... + t^{r-1}\sum\limits_{1 \leq  i_{r-1}< ... <i_2 <i_1 < r} \  \prod\limits_{\substack{j=1 \\ i_0 = r}}^{r-1}  \frac{w^{i_j}(z;q)_{i_j}(w;q)_{i_{j-1}-i_j}}{((zw;q)_{i_j}-tw^{i_j} (z;q)_{i_j})(q;q)_{i_{j-1}-i_j} } \right\}\\ = \left(\frac{1-q^n}{1-zwq^{n-1}}\right)\frac{(zw;q)_n}{(zw;q)_n-tw^n (z;q)_n} . \tag{1.13} \label{egnewexpression2}
\end{multline}

As one of the consequences, we will also arrive at the following identity for the reciprocal harmonic number : 

\begin{multline}
   \frac{1}{H_n}\left\{ 1+  \sum\limits_{1 \leq r_1 < n} \frac{1}{H_{r_1}(n-r_1)} +  \sum\limits_{1 \leq r_2 < r_1 < n} \frac{1}{H_{r_1}H_{r_2}(n-r_1)(r_1-r_2)} + ... \right. \\ \left. ...+ \sum\limits_{1 \leq r_{n-1} <...<r_2<r_1<n} \ \prod\limits_{\substack{j=1 \\ r_0 = n}}^{n-1}\frac{1}{H_{r_j}(r_{j-1}-r_j)}  \right\}  = \sum\limits_{r\geq 1} \binom{n}{r} \frac{(-1)^{r-1}}{H_r}, \tag{1.14} \label{egnewexpression3}
\end{multline}

where $H_n = \sum\limits_{i=1}^n \frac{1}{i}$ is the ordinary harmonic number.

\vspace{0.7cm}
\section{Demonstration of Theorem 1}

Let us first invoke some concepts and notations from $q$-calculus. We shall denote by $D_q$ , the $q$-derivative of a function $f$.

\[ \tag{2.1} \label{qdiff}
D_qf(x) = (D_qf)(x) = \frac{f(x)-f(xq)}{x-xq}.
\]

The definite Jackson integral of $f$ is defined in \cite{JacksonqInt} as

\[\tag{2.2} \label{qint}
\int_{0}^{x} f(t) d_qt = (1-q)\sum\limits_{n \geq 0} q^nx f(q^nx).
\]

The Jackson integral and the $q$-derivative are related by the following \emph{fundamental theorem of quantum calculus} \cite[p.~73]{KacCheung}, which implies that if $D_qF = f$ and $F$ is continuous at $x=0$, then

\[\tag{2.3} \label{qprop1}
\int_{0}^{x}f(t)d_qt = F(x)-F(0).
\]

Furthermore, for any function $f$

\[\tag{2.4} \label{qprop2}
D_q \int_{0}^{x} f(t) d_qt = f(x).
\]

Then, with all of these assumed, it is easily deduced that 

\[ \tag{2.5} \label{intx^n}
\int_{0}^{x} t^{n-1} d_qt = \frac{x^n}{[n]},
\]
and 

\[ \tag{2.6} \label{int(x)n}
\int_{0}^{x} (tyq;q)_{n-1} d_qt  = -\frac{(xy;q)_n}{y[n]}.
\]

Now, for the purpose of our demonstration, using the Jackson integral, let us further define the operators $P_q$ and $T_q$ as follows.

\[
P_q f(x) = \int_{0}^{x} \frac{f(t)}{t} d_qt, \hspace{0.05cm} P_q^m f(x) = P_q^{m-1}P_q f(x), \hspace{0.05cm}  m \in \mathbb{N}, \hspace{0.05cm} P_q^0 f(x) = f(x),
\]

and 
\[
T_q f(x) = \int_{0}^{x}\frac{1-f(t)}{1-t} d_qt, \hspace{0.05cm} T_q^{m} = T_q^{m-1} T_q f(x),\hspace{0.05cm} m \in \mathbb{N}, \hspace{0.05cm} T_q^0 f(x) = f(x).
\]

Then, we shall state the following lemma on $P_q$ and $T_q$.
\vspace{0.7cm}
\begin{lemma} For a non-negative integer $m$, the following transformations hold.

\[ \tag{2.7} \label{trans1}
P_q^m x^n = \frac{x^n}{[n]^m}, 
\]
\[ \tag{2.8} \label{trans2}
P_q^m (1-(x;q)_n) = \sum\limits_{n\geq i_1 \geq i_2 \geq ... \geq i_m \geq 1}\frac{q^{i_1+i_2+...+i_m}}{[i_1][i_2]...[i_m]}(1-(xq^{-m};q)_{i_m}), 
\]
\[ \tag{2.9} \label{trans3}
T_q^m x^n = \sum\limits_{n\geq i_1 \geq i_2 \geq ... \geq i_m \geq 1}\frac{1}{[i_1][i_2]...[i_m]}x^{i_m}, 
\]
\[ \tag{2.10} \label{trans4}
T_q^m (1-(x;q)_n) = \frac{1-(x;q)_n}{[n]^m}.
\]
    
\end{lemma}

\begin{proof}[Proof of Lemma 6]

For  \eqref{trans1}, we simply note that
\[
P_q^m x^n = P_q^{m-1} \int_{0}^{x} t^{n-1} d_qt = \frac{1}{[n]} P_q^{m-1}x^n=...= \frac{x^n}{[n]^m}. 
\]

For \eqref{trans2}, first, we note
\[
P_q^m (1-(x;q)_n) = P_q^{m-1} \int_{0}^{x} \frac{1-(t;q)_n}{t} d_qt.
\]
But since $\frac{1-(t;q)_n}{t} = \sum\limits_{n \geq i \geq 1} q^{i-1}(t;q)_{i-1}$, we have

\begin{multline} \notag
P_q^m (1-(x;q)_n) = \sum\limits_{n \geq i \geq 1} q^{i-1} P_q^{m-1} \int_{0}^{x} (t;q)_{i-1} d_qt = \sum\limits_{n \geq i \geq 1} \frac{q^i}{[i]} P_q^{m-1} (1-(xq^{-1};q)_i)\\=... = \sum\limits_{n\geq i_1 \geq i_2 \geq ... \geq i_m \geq 1}\frac{q^{i_1+i_2+...+i_m}}{[i_1][i_2]...[i_m]}(1-(xq^{-m};q)_{i_m}) .
\end{multline}

For \eqref{trans3}, we use the fact $\frac{1-t^n}{1-t} = 1+ t + ... + t^{n-1}$, to get
\[
T_q^m x^n = T_q^{m-1} \int_{0}^{x} \frac{1-t^n}{1-t} d_qt = \sum\limits_{n \geq i \geq 1} T_q^{m-1} \int_{0}^{x} t^{i-1} d_qt = \sum\limits_{n \geq i \geq 1} \frac{1}{[i]}T_q^{m-1} t^{i} =...= \sum\limits_{n\geq i_1 \geq i_2 \geq ... \geq i_m \geq 1}\frac{1}{[i_1][i_2]...[i_m]}x^{i_m}.
\]

For \eqref{trans4}, we simply evaluate as
\[
T_q^m (1-(x;q)_n) = T_q^{m-1} \int_{0}^{x} (tq;q)_{n-1} d_qt = \frac{1}{[n]} T_q^{m-1} (1-(x;q)_n) =...= \frac{1-(x;q)_n}{[n]^m}.
\]
Now, we have proved all the transformations \eqref{trans1}-\eqref{trans4}.
\end{proof}
We shall now give our proof of \textbf{Theorem 1}.

\begin{proof}[Proof of Theorem 1] Suppose that we have the sequences $A_1(q), A_2(q),...$ and $B_1(q), B_2(q),...$, not depending on $x$, satisfying the equality  
    \[
    \sum\limits_{r\geq 1} A_r(q) x^r = \sum\limits_{r \geq 1} B_r(q) (1-(x;q)_r),
    \]
for all complex values $x$. Then, let us consider applying a series of combinations of operators $(T_q^{n_k}\eta_x^{m_k}P_q^{m_k})...$\\$...(T_q^{n_2}\eta_x^{m_2}P_q^{m_2})(T_q^{n_1}\eta_x^{m_1}P_q^{m_1})$ to both sides of the above equality. In view of \textbf{Lemma 6}, after the application of the innermost combination $(T_q^{n_1}\eta_x^{m_1}P_q^{m_1})$, we arrive at 

\begin{multline} \notag
    \sum\limits_{r\geq 1} A_r(q)\frac{q^{rm_1}}{[r]^{m_1}} \sum\limits_{r\geq i_1 \geq ... \geq i_{n_1} \geq 1}\frac{1}{[i_1]...[i_{n_1}]}(T_q^{n_k}\eta_x^{m_k}P_q^{m_k})...(T_q^{n_3}\eta_x^{m_3}P_q^{m_3})(T_q^{n_2}\eta_x^{m_2}P_q^{m_2}) x^{i_{n_1}} \\ = \sum\limits_{r \geq 1} B_r(q) \sum\limits_{r\geq i_1 \geq ... \geq i_{m_1} \geq 1}\frac{q^{i_1+...+i_{m_1}}}{[i_1]...[i_{m_1}]}\frac{1}{[i_{m_1}]^{n_1}} (T_q^{n_k}\eta_x^{m_k}P_q^{m_k})...(T_q^{n_3}\eta_x^{m_3}P_q^{m_3})(T_q^{n_2}\eta_x^{m_2}P_q^{m_2})(1-(x;q)_{i_{m_1}}).
\end{multline}

Now, repeating this operation until no combinations $(T_q^{n_j}\eta_x^{m_j}P_q^{m_j})$ are left, with our definition of \emph{generalized q-multiple harmonic sums} in mind, we obtain

\begin{multline} \notag
\sum\limits_{r \geq 1} A_r(q) H_r[m_1+1, \{1\}^{n_1-1}, m_2+1, \{1\}^{n_2-1}, m_3+1,..., \{1\}^{n_{k-1}-1}, m_k+1, \{1\}^{n_k} ; x] \\= \sum\limits_{r \geq 1} B_r(q) U_r[\{1\}^{m_1}, n_1+1, \{1\}^{m_2-1}, n_2+1, \{1\}^{m_3-1}, n_3+1,..., \{1\}^{m_k-1},n_k+1;x],
\end{multline}

For all non-negative integers $m_1,m_2,...,m_k$ and $n_1,n_2,...,n_k$. And if some $n_j=0$, we will have
\begin{multline} \notag
H_r[m_1+1, \{1\}^{n_1-1}, m_2+1, \{1\}^{n_2-1}, m_3+1,...,\{1\}^{n_{j-1}-1},m_{j}+1,\{1\}^{n_{j}-1}, m_{j+1}+1,..., \{1\}^{n_{k-1}-1}, m_k+1, \{1\}^{n_k} ; x]\\ = H_r[m_1+1, \{1\}^{n_1-1}, m_2+1, \{1\}^{n_2-1}, m_3+1,...,\{1\}^{n_{j-1}-1},m_{j}+m_{j+1}+1,..., \{1\}^{n_{k-1}-1}, m_k+1, \{1\}^{n_k} ; x],
\end{multline}

since when $n_j=0$

\begin{multline} \notag
(T_q^{n_k}\eta_x^{m_k}P_q^{m_k})...(T_q^{n_{j+1}}\eta_x^{m_{j+1}}P_q^{m_{j+1}})(T_q^{n_{j}}\eta_x^{m_{j}}P_q^{m_{j}})(T_q^{n_{j-1}}\eta_x^{m_{j-1}}P_q^{m_{j-1}})...(T_q^{n_1}\eta_x^{m_1}P_q^{m_1})x^r \\= (T_q^{n_k}\eta_x^{m_k}P_q^{m_k})...(T_q^{n_{j+2}}\eta_x^{m_{j+2}}P_q^{m_{j+2}})(T_q^{n_{j+1}}\eta_x^{m_{j+1}+m_j}P_q^{m_{j+1}+m_j})(T_q^{n_{j-1}}\eta_x^{m_{j-1}}P_q^{m_{j-1}})...(T_q^{n_1}\eta_x^{m_1}P_q^{m_1})x^r.
\end{multline}

We can observe that the same holds for $U_r$. Thus, we have completed our proof of \textbf{Theorem 1}.

\end{proof}

The most obvious application of \textbf{Theorem 1} would be the application to the $q$-binomial theorem in the following form. 

\[ \tag{2.11} \label{gaussbinom}
\sum\limits_{r\geq 1} \gaussian{n}{r} (-1)^{r-1}q^{\binom{r}{2}}x^r = 1- (x;q)_n.
\]

Then we have $A_r(q) = \gaussian{n}{r} (-1)^{r-1}q^{\binom{r}{2}}$ and $B_n(q) =1$, $B_r(q) = 0$ for all $r \neq n$. Thus, we state the following corollary. 
\vspace{0.7cm}
\begin{corollary} For all non-negative integers $m_1,m_2,...,m_k$ and $n_1,n_2,...,n_k$, there holds
    \begin{multline} \tag{2.12} \label{multigauss}
\sum\limits_{r \geq 1}  \gaussian{n}{r} (-1)^{r-1}q^{\binom{r}{2}} H_r[m_1+1, \{1\}^{n_1-1}, m_2+1, \{1\}^{n_2-1}, m_3+1,..., \{1\}^{n_{k-1}-1}, m_k+1, \{1\}^{n_k} ; x] \\=  U_n[\{1\}^{m_1}, n_1+1, \{1\}^{m_2-1}, n_2+1, \{1\}^{m_3-1}, n_3+1,..., \{1\}^{m_k-1},n_k+1;x].
\end{multline}
Where it is understood that when $s=0$, $H_N[n_1,...,n_p, n+1, \{1\}^{s-1}, r+1,r_1,...,r_q;x] = H_N[n_1,...,n_p, n+r+1,r_1,...,r_q;x]$ and $U_N[n_1,...,n_p, n+1, \{1\}^{s-1}, r+1,r_1,...,r_q;x] = U_N[n_1,...,n_p, n+r+1,r_1,...,r_q;x]$. 
\end{corollary}

Note that \textbf{Corollary 7} is essentially equivalent to Theorem A given by the author in \cite{maw2024qseriesidentitiesmultiplesums}. When $x=1$ \textbf{Corollary 7} gives the duality relation (Theorem 1) given by Bradley \cite{Bradley}. Let $k=1$ to get 

\[ \tag{2.13} \label{multigausscase1}
\sum\limits_{r \geq 1}  \gaussian{n}{r} (-1)^{r-1}q^{\binom{r}{2}}  H_r[m_1+1, \{1\}^{n_1};x] = U_n[\{1\}^{m_1}, n_1+1;x].
\]

When $m_1=m$, $n_1=0$, we get 

\[\tag{2.14} \label{multigausscase2}
\sum\limits_{r \geq 1}  \gaussian{n}{r} \frac{(-1)^{r-1}q^{\binom{r}{2}+rm}}{[r]^m}x^r = \sum\limits_{n\geq i_1 \geq ... \geq i_m \geq 1} (1-(x;q)_{i_m})\frac{q^{i_1+...+i_m}}{[i_1]...[i_m]}.
\]

When \eqref{multigausscase2} is transformed using \textbf{Lemma 4}, we see that it is equivalent to the result of Fu and Lascoux \eqref{FuLascoux}. If we put $n_1 = m$ and $m_1 = 1$, we get 

\[\tag{2.15} \label{multigausscase3}
\sum\limits_{r \geq 1}  \gaussian{n}{r} \frac{(-1)^{r-1}q^{\binom{r+1}{2}}}{[r]}\sum\limits_{r\geq i_1 \geq ... \geq i_m \geq 1} \frac{x^{i_m}}{[i_1]...[i_m]}= \sum\limits_{n\geq r \geq 1} (1-(x;q)_{r})\frac{q^{r}}{[r]^{m+1}}.
\]

Which generalizes the result of Prodinger \eqref{prodinger}. We proceed to give a few further examples. In \eqref{multigausscase1}, put $m_1=2$ and $n_1=3$ to get

\[ \tag{2.16} \label{multigausscase4}
\sum\limits_{r \geq 1}  \gaussian{n}{r} \frac{(-1)^{r-1}q^{\binom{r+1}{2}+r}}{[r]^2} \sum\limits_{r \geq i_1 \geq i_2 \geq i_3 \geq 1} \frac{x^{i_3}}{[i_1][i_2][i_3]} = \sum\limits_{n \geq i_1 \geq i_2 \geq 1} (1-(x;q)_{i_2})\frac{q^{i_1+i_2}}{[i_1][i_2]^4}, 
\]

and when $m_1 = 3$, $n_1=2$, we get

\[ \tag{2.17} \label{multigausscase5}
\sum\limits_{r \geq 1}  \gaussian{n}{r} \frac{(-1)^{r-1}q^{\binom{r+1}{2}+2r}}{[r]^3} \sum\limits_{r \geq i_1 \geq i_2 \geq 1} \frac{x^{i_2}}{[i_1][i_2]} = \sum\limits_{n \geq i_1 \geq i_2 \geq i_3\geq 1} (1-(x;q)_{i_3})\frac{q^{i_1+i_2+i_3}}{[i_1][i_2][i_3]^3}. 
\]
As the last example, we put $k=2$, $m_1=2$, $n_1=3$, $m_2=2$, $n_2=1$, in \textbf{Corollary 7} to get

\[ \tag{2.18} \label{multigausscase6}
\sum\limits_{r \geq 1}  \gaussian{n}{r} \frac{(-1)^{r-1}q^{\binom{r+1}{2}+r}}{[r]^2} \sum\limits_{r \geq i_1 \geq i_2 \geq i_3 \geq i_4 \geq 1} \frac{x^{i_4}q^{2i_3}}{[i_1][i_2][i_3]^3[i_4]} = \sum\limits_{n \geq i_1 \geq i_2 \geq i_3 \geq i_4 \geq 1} (1-(x;q)_{i_4})\frac{q^{i_1+i_2+i_3+i_4}}{[i_1][i_2]^4[i_3][i_4]^2}. 
\]
\vspace{0.7cm}
\section{Demonstration of Theorem 2}
Let us first recall the following definition of the $k^{th}$ complete symmetric function $h_k$,

\[
h_k(a_1,...,a_n) = \sum\limits_{n \geq i_1 \geq ... \geq i_k \geq 1} a_{i_1}...a_{i_k}, \hspace{0.05cm} \text{with } h_0(a_1,...,a_n) =1.
\]

Then, the generating function of $h_k$ is given as

\[ \tag{3.1} \label{symmetriceq}
\sum\limits_{k \geq 0} t^k h_k(a_1,...,a_n) = \frac{1}{(1-a_1t)...(1-a_nt)}. 
\]

Now, for our purpose, we shall state the following lemma.
\vspace{0.7cm}
\begin{lemma} For an arbitrary sequence $a_1,...,a_n$, and for a complex value $z$, the following transformations hold.

\[ \tag{3.2} \label{symmetricetrans1}
\sum\limits_{k \geq 1} \frac{(z-1)^{k-1}}{(1-q)^k} \sum\limits_{n \geq i_1 \geq ... \geq i_k \geq 1} \frac{q^{i_1+...+i_k}}{[i_1]...[i_k]} a_{i_k} = \frac{(q;q)_n}{(zq;q)_n} \sum\limits_{n \geq i \geq 1} \frac{q^i(zq;q)_{i-1}}{(q;q)_i}a_i.
\]

\[ \tag{3.3} \label{symmetricetrans2}
\sum\limits_{k \geq 1} \frac{z^{-1}(1-z^{-1})^{k-1}}{(1-q)^k} \sum\limits_{n \geq i_1 \geq ... \geq i_k \geq 1} \frac{1}{[i_1]...[i_k]} a_{i_k} = \frac{(q;q)_n}{(zq;q)_n} \sum\limits_{n \geq i \geq 1} \frac{z^{n-i}(zq;q)_{i-1}}{(q;q)_i}a_i.
\]

\[\tag{3.4} \label{symmetricetrans3}
\sum\limits_{k \geq 1} \frac{(z-1)^{k-1}}{(1-q)^k} \frac{q^{ik}}{[i]^k} = \frac{q^i}{1-zq^i}.
\]

\[\tag{3.5} \label{symmetricetrans4}
\sum\limits_{k \geq 1} \frac{z^{-1}(1-z^{-1})^{k-1}}{(1-q)^k} \frac{1}{[i]^k} = \frac{1}{1-zq^i}.
\]

Provided that all the expressions from both left-hand and right-hand sides converges.
\end{lemma}

\begin{proof}[Proof of Lemma 8] Now in \eqref{symmetriceq}, if we replace the sequence $a_1,...,a_n$ with the sequence $\frac{q^i}{1-q^i},..., \frac{q^n}{1-q^n}$, and put $t=z-1$, we arrive to the following fact

\[\tag{3.6} \label{principlesymmetrictrans}
\sum\limits_{k \geq 1} (z-1)^{k-1} h_{k-1}\left(\frac{q^i}{1-q^i},...,\frac{q^n}{1-q^n} \right) = \frac{1-q^i}{1-zq^i}\frac{(q;q)_n (zq;q)_i}{(zq;q)_n (q;q)_i}
\]
To prove \eqref{symmetricetrans1}, we multiply both sides of the above equality by $\frac{q^i}{1-q^i}a_i$ and sum through $i=1,2,..,n$, to get

\[
\sum\limits_{k \geq 1} (z-1)^{k-1} \sum\limits_{n \geq i \geq 1} a_i\frac{q^i}{1-q^i} h_{k-1}\left(\frac{q^i}{1-q^i},...,\frac{q^n}{1-q^n}\right) = \frac{(q;q)_n}{(zq;q)_n} \sum\limits_{n \geq i \geq 1} \frac{q^i(zq;q)_{i-1}}{(q;q)_i}a_i.
\]
But since 
\[ \sum\limits_{n \geq i \geq 1} a_i\frac{q^i}{1-q^i} h_{k-1}\left(\frac{q^i}{1-q^i},...,\frac{q^n}{1-q^n}\right) = \sum\limits_{n \geq i_1 \geq ... \geq i_k \geq 1} \frac{q^{i_1+...+i_k}}{(1-q^{i_1})...(1-q^{i_k})} a_{i_k}, \]

we arrive at \eqref{symmetricetrans1}. To arrive at $(3.3)$, we replace $z$ by $z^{-1}$ and $q$ by $q^{-1}$ in \eqref{symmetricetrans1}. \eqref{symmetricetrans3} and \eqref{symmetricetrans4} are deduced by summing the series from the left-hand side as geometric series.
\end{proof}

We shall now demonstrate our proof of \textbf{Theorem 2}.

\begin{proof}[Proof of Theorem 2] 
Suppose that we have the sequences $A_1(q), A_2(q),...$ and $B_1(q), B_2(q),...$, not depending on $x$, satisfying the equality  
    \[
    \sum\limits_{r\geq 1} A_r(q) x^r = \sum\limits_{r \geq 1} B_r(q) (1-(x;q)_r),
    \]
for all complex values $x$. Then by \textbf{Theorem 1}, we recall that we have the following equality 
\begin{multline} \notag
\sum\limits_{r \geq 1} A_r(q) H_r[m_1+1, \{1\}^{n_1-1}, m_2+1, \{1\}^{n_2-1}, m_3+1,..., \{1\}^{n_{k-1}-1}, m_k+1, \{1\}^{n_k} ; x] \\= \sum\limits_{r \geq 1} B_r(q) U_r[\{1\}^{m_1}, n_1+1, \{1\}^{m_2-1}, n_2+1, \{1\}^{m_3-1}, n_3+1,..., \{1\}^{m_k-1},n_k+1;x].
\end{multline}

Then we also have
\begin{multline}\notag
    \sum\limits_{r \geq 1} A_r(q)\sum\limits_{\substack{m_j \geq 1 \\ 1 \leq j \leq k}}\frac{(z_1-1)^{m_1-1}...(z_k-1)^{m_k-1}}{(1-q)^{m_1+...+m_k}}\sum\limits_{\substack{n_l \geq 1 \\ 1 \leq l \leq k}}\frac{y_1^{-1}(1-y_1^{-1})^{n_1-1}...y_k^{-1}(1-y_k^{-1})^{n_k-1}}{(1-q)^{n_1+...+n_k}} \\H_r[m_1+1, \{1\}^{n_1-1}, m_2+1, \{1\}^{n_2-1}, m_3+1,..., \{1\}^{n_{k-1}-1}, m_k+1, \{1\}^{n_k} ; x] \\= \sum\limits_{r \geq 1} B_r(q)\sum\limits_{\substack{m_j \geq 1 \\ 1 \leq j \leq k}}\frac{(z_1-1)^{m_1-1}...(z_k-1)^{m_k-1}}{(1-q)^{m_1+...+m_k}}\sum\limits_{\substack{n_l \geq 1 \\ 1 \leq l \leq k}}\frac{y_1^{-1}(1-y_1^{-1})^{n_1-1}...y_k^{-1}(1-y_k^{-1})^{n_k-1}}{(1-q)^{n_1+...+n_k}} \\ U_r[\{1\}^{m_1}, n_1+1, \{1\}^{m_2-1}, n_2+1, \{1\}^{m_3-1}, n_3+1,..., \{1\}^{m_k-1},n_k+1;x].
\end{multline}

Summing over the outermost pair $m_1$ and $n_1$, with \textbf{Lemma 8} in mind, we have

\begin{multline} \notag
    \sum\limits_{r \geq 1} A_r(q)\frac{q^r(q;q)_r}{(1-z_1q^r)(y_1q;q)_r}\sum\limits_{r \geq r_1 \geq 1} \frac{y_1^{r-r_1}(y_1q;q)_{r_1-1}}{(q;q)_{r_1}} \sum\limits_{\substack{m_j \geq 1 \\ 2 \leq j \leq k}}\frac{(z_2-1)^{m_2-1}...(z_k-1)^{m_k-1}}{(1-q)^{m_2+...+m_k}}\\\sum\limits_{\substack{n_l \geq 1 \\ 2 \leq l \leq k}}\frac{y_2^{-1}(1-y_2^{-1})^{n_2-1}...y_k^{-1}(1-y_k^{-1})^{n_k-1}}{(1-q)^{n_2+...+n_k}} H_{r_1}[ m_2+1, \{1\}^{n_2-1}, m_3+1,..., \{1\}^{n_{k-1}-1}, m_k+1, \{1\}^{n_k} ; x] \\= \sum\limits_{r \geq 1} B_r(q) \frac{(q;q)_r}{(z_1q;q)_r} \sum\limits_{r \geq r_1 \geq 1} \frac{q^{r_1}(z_1q;q)_{r_1-1}}{(1-y_1q^{r_1})(q;q)_{r_1}} \sum\limits_{\substack{m_j \geq 1 \\ 2 \leq j \leq k}}\frac{(z_2-1)^{m_2-1}...(z_k-1)^{m_k-1}}{(1-q)^{m_2+...+m_k}}\\\sum\limits_{\substack{n_l \geq 1 \\ 2 \leq l \leq k}}\frac{y_2^{-1}(1-y_2^{-1})^{n_2-1}...y_k^{-1}(1-y_k^{-1})^{n_k-1}}{(1-q)^{n_2+...+n_k}} U_{r_1}[\{1\}^{m_2}, n_2+1, \{1\}^{m_3-1}, n_3+1,..., \{1\}^{m_k-1},n_k+1;x].
\end{multline}

Now, repeating this operation until no pairs $m_j$ and $n_j$ are left, we arrive at 
\begin{multline} \notag
    \sum\limits_{r\geq 1} A_r(q)\frac{q^r(q;q)_r}{(1-z_1q^r)(y_1q;q)_r}\sum\limits_{r \geq r_1 \geq 1}\frac{q^{r_1}y_1^{r-r_1}(y_1q;q)_{r_1-1}}{(1-z_2q^{r_1})(y_2q;q)_{r_1}}\sum\limits_{r_1 \geq r_2 \geq 1} \frac{q^{r_2}y_2^{r_1-r_2}(y_2q;q)_{r_2-1}}{(1-z_3q^{r_2})(y_3q;q)_{r_2}}\hspace{0.05cm}...\\...\hspace{0.05cm}\sum\limits_{r_{k-2} \geq r_{k-1} \geq 1} \frac{q^{r_{k-1}}y_{k-1}^{r_{k-2}-r_{k-1}} (y_{k-1}q;q)_{r_{k-1}-1}}{(1-z_{k}q^{r_{k-1}})(y_kq;q)_{r_{k-1}}}\sum\limits_{r_{k-1} \geq r_{k} \geq 1} \frac{x^{r_{k}}y_{k}^{r_{k-1}-r_{k}}(y_kq;q)_{r_k-1}}{(q;q)_{r_k}} 
\end{multline}

\begin{multline} \notag
= \sum\limits_{r\geq 1} B_r(q) \frac{(q;q)_r}{(z_1q;q)_r}\sum\limits_{r \geq r_1 \geq 1}\frac{q^{r_1}(z_1q;q)_{r_1-1}}{(1-y_1q^{r_1})(z_2q;q)_{r_1}}\sum\limits_{r_1 \geq r_2 \geq 1}\frac{q^{r_2}(z_2q;q)_{r_2-1}}{(1-y_2q^{r_2})(z_3q;q)_{r_2}}\hspace{0.05cm}...\\...\hspace{0.05cm}\sum\limits_{r_{k-2} \geq r_{k-1} \geq 1}\frac{q^{r_{k-1}}(z_{k-1}q;q)_{r_{k-1}-1}}{(1-y_{k-1}q^{r_{k-1}})(z_kq;q)_{r_{k-1}}}\sum\limits_{r_{k-1} \geq r_{k} \geq 1}\frac{(1-(x;q)_{r_k})q^{r_{k}}(z_{k}q;q)_{r_{k}-1}}{(1-y_{k}q^{r_{k}})(q;q)_{r_{k}}}.
\end{multline}

\end{proof}

The application of \textbf{Theorem 2} to the $q$-binomial theorem in the form \eqref{gaussbinom}, with $A_r(q) = \gaussian{n}{r} (-1)^{r-1}q^{\binom{r}{2}}$ and with $B_n(q) =1$, and $B_r(q) = 0$ for all $r \neq n$, gives the following corollary.
\vspace{0.7cm}
\begin{corollary} For complex values $x$, $z_1,z_2,...,z_k$ and $y_1,y_2,...,y_k$, there holds

\begin{multline} \notag
    \sum\limits_{r\geq 1} \gaussian{n}{r} \frac{(-1)^{r-1}q^{\binom{r+1}{2}}(q;q)_r}{(1-z_1q^r)(y_1q;q)_r}\sum\limits_{r_0=r \geq r_1 \geq r_2 \geq \cdots \geq r_{k-1} \geq r_k \geq 1} \frac{x^{r_{k}}y_{k}^{r_{k-1}-r_{k}}(y_kq;q)_{r_k-1}}{(q;q)_{r_k}} \prod\limits_{j=1}^{k-1} \frac{q^{r_j}y_j^{r_{j-1}-r_j}(y_jq;q)_{r_j-1}}{(1-z_{j+1}q^{r_j})(y_{j+1}q;q)_{r_j}} 
\end{multline}

\begin{multline} \tag{3.7} \label{analogguassainmultiple}
= \frac{(q;q)_n}{(z_1q;q)_n}\sum\limits_{n \geq r_1 \geq r_2 \geq \cdots \geq r_{k-1} \geq r_k \geq 1} \frac{(1-(x;q)_{r_k})q^{r_{k}}(z_{k}q;q)_{r_{k}-1}}{(1-y_{k}q^{r_{k}})(q;q)_{r_{k}}} \prod\limits_{j=1}^{k-1} \frac{q^{r_j}(z_jq;q)_{r_j-1}}{(1-y_jq^{r_j})(z_{j+1}q;q)_{r_j}}.
\end{multline}
Provided that $z_j, y_j \neq q^{-r}$, $r \in \mathbb{N}$.
\end{corollary}

When we let $n \to \infty$, this gives the following $q$-multiple infinite sum identity.
\vspace{0.7cm}
\begin{corollary} For complex values $|q| < 1$, $x$, $z_1,z_2,...,z_k$ and $y_1,y_2,...,y_k$, there holds
    \begin{multline} \notag
    \sum\limits_{r=1}^{\infty} \frac{(-1)^{r-1}q^{\binom{r+1}{2}}}{(1-z_1q^r)(y_1q;q)_r}\sum\limits_{r_0=r \geq r_1 \geq r_2 \geq \cdots \geq r_{k-1} \geq r_k \geq 1} \frac{x^{r_{k}}y_{k}^{r_{k-1}-r_{k}}(y_kq;q)_{r_k-1}}{(q;q)_{r_k}} \prod\limits_{j=1}^{k-1} \frac{q^{r_j}y_j^{r_{j-1}-r_j}(y_jq;q)_{r_j-1}}{(1-z_{j+1}q^{r_j})(y_{j+1}q;q)_{r_j}} 
\end{multline}

\begin{multline} \tag{3.8} \label{infiniteanalogguassainmultiple}
= \frac{(q;q)_\infty}{(z_1q;q)_\infty}\sum\limits_{ r_1 \geq r_2 \geq \cdots \geq r_{k-1} \geq r_k \geq 1} \frac{(1-(x;q)_{r_k})q^{r_{k}}(z_{k}q;q)_{r_{k}-1}}{(1-y_{k}q^{r_{k}})(q;q)_{r_{k}}} \prod\limits_{j=1}^{k-1} \frac{q^{r_j}(z_jq;q)_{r_j-1}}{(1-y_jq^{r_j})(z_{j+1}q;q)_{r_j}}.
\end{multline}
Provided that $z_j, y_j \neq q^{-r}$, $r \in \mathbb{N}$.
\end{corollary}

Now let $z_1=...=z_k=z$ and $y_1=...=y_k=y$ in \eqref{analogguassainmultiple} to get 

\vspace{0.7cm}
\begin{corollary} For complex values $x$, $y$ and $z$, there holds

\begin{multline} \notag
    \sum\limits_{r\geq 1} \gaussian{n}{r} \frac{(-1)^{r-1}y^rq^{\binom{r+1}{2}}(q;q)_r}{(1-zq^r)(yq;q)_r}\sum\limits_{r \geq r_1 \geq r_2 \geq \cdots \geq r_{k-1} \geq r_k \geq 1} \frac{x^{r_{k}}y^{-r_{k}}(yq;q)_{r_k-1}}{(q;q)_{r_k}} \prod\limits_{j=1}^{k-1}\frac{q^{r_j}}{(1-zq^{r_j})(1-yq^{r_j})}  
\end{multline}

\begin{multline} \tag{3.9} \label{analogguassainmultiple2}
= \frac{(q;q)_n}{(zq;q)_n}\sum\limits_{n \geq r_1 \geq r_2 \geq \cdots \geq r_{k-1} \geq r_k \geq 1}\frac{(1-(x;q)_{r_k})q^{r_{k}}(zq;q)_{r_{k}-1}}{(1-yq^{r_{k}})(q;q)_{r_{k}}} \prod\limits_{j=1}^{k-1} \frac{q^{r_j}}{(1-zq^{r_j})(1-yq^{r_j})}.
\end{multline}

Provided that $y, z \neq q^{-r}$, $r \in \mathbb{N}$.
\end{corollary}

When, $y=0$, this reduces into

\begin{multline} \tag{3.10} \label{analoggaussianmultiple3}
    \sum\limits_{r\geq 1} \gaussian{n}{r} \frac{(-1)^{r-1}x^rq^{\binom{r}{2}+rk}}{(1-zq^r)^{k}} = \frac{(q;q)_n}{(zq;q)_n}\sum\limits_{n \geq r_1 \geq r_2 \geq \cdots \geq r_{k-1} \geq r_k \geq 1}\frac{(1-(x;q)_{r_k})q^{r_{k}}(zq;q)_{r_{k}-1}}{(q;q)_{r_{k}}} \prod\limits_{j=1}^{k-1} \frac{q^{r_j}}{1-zq^{r_j}}.
\end{multline}

Replacing $z$ by $z^{-1}$ and using \textbf{Lemma 4} this becomes

\begin{multline} \tag{3.11} \label{analoggaussianmultiple3alt}
    \sum\limits_{r\geq 1} \gaussian{n}{r} \frac{(-1)^{r-1}(1-(x;q)_r)q^{\binom{r+1}{2}-rn}}{(1-zq^r)^{k}} = \frac{(q;q)_n}{(zq;q)_n}\sum\limits_{n \geq r_1 \geq r_2 \geq \cdots \geq r_{k-1} \geq r_{k} \geq 1} \frac{x^{r_k}z^{n-r_k}(zq;q)_{r_{k}-1}}{(q;q)_{r_{k}}} \prod\limits_{j=1}^{k-1} \frac{1}{1-zq^{r_j}}.
\end{multline}

This is equivalent to \eqref{Zeng} by Zeng. When we apply \textbf{Lemma 5} to \eqref{analoggaussianmultiple3}, we arrive at

\begin{multline} \tag{3.12} \label{analoggaussianmultiple3extension}
    \sum\limits_{r\geq 1} \gaussian{n}{r} \left(1-x^r\frac{(y;q)_r}{(xy;q)_r}\right)\frac{(-1)^{r-1}q^{\binom{r}{2}+rk}}{(1-zq^r)^{k}} = \frac{(q;q)_n}{(zq;q)_{n}}\sum\limits_{n \geq r_1 \geq r_2 \geq \cdots \geq r_{k-1} \geq r_k \geq 1} \frac{q^{r_{k}}(x;q)_{r_k}(zq;q)_{r_{k}-1}}{(q;q)_{r_{k}}(xy;q)_{r_k}} \prod\limits_{j=1}^{k-1} \frac{q^{r_j}}{1-zq^{r_j}}.
\end{multline}
\vspace{0.7cm}
\section{Further \texorpdfstring{$q$}{q}-combinatorial results}

Based on our established results, we proceed to state the following lemma which the author was unable to find in $q$-literature.
\vspace{0.7cm}
\begin{lemma}
    For natural number $n$, and for complex values $x,w,y,z$, there holds
    \[ \tag{4.1} \label{binomtransform}
    \frac{1}{(yz;q)_n}\sum\limits_{r \geq 1} \gaussian{n}{r}\frac{x^rz^{n-r}(w;q)_r(z;q)_r(y;q)_{n-r}}{(xw;q)_r} = \sum\limits_{r \geq 1} \gaussian{n}{r}(-1)^{r-1}q^{\binom{r+1}{2}-rn}\left(1-\frac{(x;q)_r}{(xw;q)_r}\right)\frac{(z;q)_r}{(yz;q)_r}.
    \]
    Provided that $xw, yz \neq q^{-r}$, $r \in \mathbb{N} \cup \{0\}$.
\end{lemma}

\begin{proof}[Proof of Lemma 12] Let $k=1$ in \eqref{analoggaussianmultiple3alt} to arrive at
    \[ \tag{4.2} \label{FuLascoux2}
    \sum\limits_{r\geq 1} \gaussian{n}{r} \frac{(-1)^{r-1}(1-(x;q)_r)q^{\binom{r+1}{2}-rn}}{1-zq^r} = \frac{(q;q)_n}{(zq;q)_n}\sum\limits_{n \geq r \geq 1}\frac{x^{r}z^{n-r}(zq;q)_{r-1}}{(q;q)_{r}}.
\]
We note that this is an identity due to Fu and Lascoux \cite{FuLascoux}. Now apply the operator $\sum\limits_{j \geq 0} \frac{(-1)^j q^{\binom{j+1}{2}}}{(q;q)_j}\eta_z^j$ to both sides of the equality to arrive at

\begin{multline} \notag
\sum\limits_{r\geq 1} \gaussian{n}{r} (-1)^{r-1}(1-(x;q)_r)q^{\binom{r+1}{2}-rn} \sum\limits_{j \geq 0} \frac{(-1)^j q^{\binom{j+1}{2}}}{(q;q)_j(1-zq^{r+j})}\\=\frac{(q;q)_n}{(zq;q)_n}\sum\limits_{n \geq r \geq 1}\frac{x^{r}z^{n-r}(zq;q)_{r-1}}{(q;q)_{r}}\sum\limits_{j \geq 0}\frac{(-1)^j q^{\binom{j+1}{2}}q^{j(n-r)}(zq^r;q)_j}{(q;q)_j(zq^{n+1})_j}.
\end{multline}

Now since \[ \sum\limits_{j \geq 0} \frac{(-1)^jq^{\binom{j+1}{2}}}{(q;q)_j(1-zq^{j})}=\frac{(q;q)_\infty}{(z;q)_\infty}, \] and \[ \sum\limits_{j \geq 0} \frac{(-1)^jq^{\binom{j+1}{2}}q^{j(n-r)}(zq^r;q)_j}{(q;q)_j (zq^{n+1};q)_j}=\frac{(q^{n-r+1};q)_\infty}{(zq^{n+1};q)_\infty}.\] For the two equalities above, we note that the former can be deduced from the partial fraction decomposition of $\frac{1}{(z;q)_\infty}$, while the later is the case $t \to 1, b = zq^n, a = q^{n-i}$ of the identity (12.2) appearing in \cite[p.~13]{nathanJ}. Therefore, we arrive at

 \[ \tag{4.3} \label{lesserbinomtransform}
\sum\limits_{r\geq 1} \gaussian{n}{r} (-1)^{r-1}(1-(x;q)_r)q^{\binom{r+1}{2}-rn} (z;q)_r = \sum\limits_{ r \geq 1} \gaussian{n}{r}x^{r}z^{n-r}(z;q)_{r} .
\]

Then, apply the operator $\sum\limits_{j \geq 0} \frac{(z;q)_j}{(q;q)_j}y^j\eta_z^j$ to both sides of the equality and use Heine's $q$-binomial theorem to get

\[ \tag{4.4} \label{lesserbinomtransform2}
\sum\limits_{r\geq 1} \gaussian{n}{r} (-1)^{r-1}(1-(x;q)_r)q^{\binom{r+1}{2}-rn} \frac{(z;q)_r}{(yz;q)_r} = \frac{1}{(yz;q)_n}\sum\limits_{ r \geq 1} \gaussian{n}{r}x^{r}z^{n-r}(z;q)_{r}(y;q)_{n-r} .
\]

Next, we use \textbf{Lemma 5} with the introduction of the new variable $w$, to finally arrive at \eqref{binomtransform}.

\end{proof}

Now, in \eqref{lesserbinomtransform2}, expanding the left-hand side in powers of $x$ and equating the coefficients of $x^i$, we arrive at the following fact.
\vspace{0.7cm}
\begin{lemma}
    For natural numbers $n$, $i$, with $i \leq n$, and for complex values $y, z$, there holds
    \[ \tag{4.5} \label{lesserbinomtransform3}
\sum\limits_{r\geq 1} \gaussian{n}{r} \gaussian{r}{i} (-1)^{r-1}q^{\binom{r+1}{2}-rn} \frac{(z;q)_r}{(yz;q)_r} = (-1)^{i-1}z^{n-i}q^{-\binom{i}{2}}\gaussian{n}{i}\frac{(z;q)_{i}(y;q)_{n-i}}{(yz;q)_n}.
\]
Provided that $yz \neq q^{-r}$, $r \in \mathbb{N} \cup \{0\}$.
\end{lemma}

Now we shall state our proof of \textbf{Proposition 3}.

\begin{proof}[Proof of Proposition 3]
Let us invoke the following inverted form of the $q$-binomial theorem

\[ \tag{4.6} \label{invertgaussbinom}
\sum\limits_{r\geq 1} \gaussian{n}{r} (-1)^{r-1}q^{\binom{r+1}{2}-rn}(1-(x;q)_r) = x^n.
\]

We shall apply \textbf{Theorem 2} to \eqref{invertgaussbinom} with $B_r(q)=\gaussian{n}{r} (-1)^{r-1}q^{\binom{r+1}{2}-rn}$ and $A_n(q) = 1$, $A_r(q)=0$ for all $r \neq n$. Then, we have 

\begin{multline} \notag
 \frac{q^n(q;q)_n}{(1-z_1q^n)(y_1q;q)_n}\sum\limits_{r_0 = n \geq r_1 \geq r_2 \geq \cdots \geq r_{k-1} \geq r_k \geq 1} \frac{x^{r_{k}}y_{k}^{r_{k-1}-r_{k}}(y_kq;q)_{r_k-1}}{(q;q)_{r_k}} \prod\limits_{j=1}^{k-1}\frac{q^{r_j}y_j^{r_{j-1}-r_j}(y_jq;q)_{r_j-1}}{(1-z_{j+1}q^{r_j})(y_{j+1}q;q)_{r_j}}
\end{multline}

\begin{multline} \notag
= \sum\limits_{r\geq 1} \gaussian{n}{r} (-1)^{r-1}q^{\binom{r+1}{2}-rn} \frac{(q;q)_r}{(z_1q;q)_r}\sum\limits_{r_0 =r \geq r_1 \geq r_2 \geq \cdots \geq r_{k-1} \geq r_k \geq 1} \frac{(1-(x;q)_{r_k})q^{r_{k}}(z_{k}q;q)_{r_{k}-1}}{(1-y_{k}q^{r_{k}})(q;q)_{r_{k}}} \prod\limits_{j=1}^{k-1} \frac{q^{r_j}(z_jq;q)_{r_j-1}}{(1-y_jq^{r_j})(z_{j+1}q;q)_{r_j}}.
\end{multline}

    Let us put $z_j = z q^{-j+1}$, and $y_j = y$, for all $1\leq j \leq k$, and make a few manipulations to arrive at

\begin{multline} \notag
 \frac{q^ny^n(q;q)_n}{(1-zq^n)(yq;q)_n}\sum\limits_{n \geq r_1 \geq r_2 \geq \cdots \geq r_{k-1} \geq r_k \geq 1} \frac{(1-x^{r_{k}})y^{-r_{k}}(yq;q)_{r_k-1}}{(q;q)_{r_k}} \prod\limits_{j=1}^{k-1} \frac{q^{r_j}}{(1-zq^{r_j-j})(1-yq^{r_j})}
\end{multline} 
\begin{multline} \notag
= \frac{1}{(zq^{-k+1};q)_k}\sum\limits_{r\geq 1} \gaussian{n}{r} (-1)^{r-1}q^{\binom{r+1}{2}-rn} \frac{(q;q)_r}{(zq;q)_r}\sum\limits_{r \geq r_1 \geq r_2 \geq \cdots \geq r_{k-1} \geq r_k \geq 1} \frac{(x;q)_{r_k}q^{r_{k}}(zq^{-k+1};q)_{r_{k}}}{(1-yq^{r_{k}})(q;q)_{r_{k}}}\prod\limits_{j=1}^{k-1}\frac{q^{r_j}}{1-yq^{r_j}}.
\end{multline}
Now let us transform this equality with the application of \textbf{Lemma 5}, introducing a new parameter $t$. Then, we have

\begin{multline} \notag
 \frac{q^ny^n(q;q)_n}{(1-zq^n)(yq;q)_n}\sum\limits_{n \geq r_1 \geq r_2 \geq \cdots \geq r_{k-1} \geq r_k \geq 1} \left(1-x^{r_{k}}\frac{(t;q)_{r_k}}{(xt;q)_{r_k}}\right)\frac{y^{-r_{k}}(yq;q)_{r_k-1}}{(q;q)_{r_k}} \prod\limits_{j=1}^{k-1}\frac{q^{r_j}}{(1-zq^{r_j-j})(1-yq^{r_j})} 
\end{multline}

\begin{multline} \notag
= \frac{1}{(zq^{-k+1};q)_k}\sum\limits_{r\geq 1} \gaussian{n}{r} (-1)^{r-1}q^{\binom{r+1}{2}-rn} \frac{(q;q)_r}{(zq;q)_r}\sum\limits_{r \geq r_1 \geq r_2 \geq \cdots \geq r_{k-1} \geq r_k \geq 1} \frac{q^{r_{k}}(x;q)_{r_k}(zq^{-k+1};q)_{r_{k}}}{(1-yq^{r_{k}})(xt;q)_{r_k}(q;q)_{r_{k}}} \prod\limits_{j=1}^{k-1} \frac{q^{r_j}}{1-yq^{r_j}}.
\end{multline}

Now let $x=yq$ and use \eqref{analoggaussianmultiple3extension} to arrive at

\begin{multline} \notag
 \frac{q^ny^n(q;q)_n}{(1-zq^n)(yq;q)_n}\sum\limits_{n \geq r_1 \geq r_2 \geq \cdots \geq r_{k-1} \geq r_k \geq 1} \left(1-y^{r_{k}}q^{r_k}\frac{(t;q)_{r_k}}{(ytq;q)_{r_k}}\right)\frac{y^{-r_{k}}(yq;q)_{r_k-1}}{(q;q)_{r_k}} \prod\limits_{j=1}^{k-1} \frac{q^{r_j}}{(1-zq^{r_j-j})(1-yq^{r_j})}  
\end{multline}

\begin{multline} \notag
= \frac{1}{(zq^{-k+1};q)_k}\sum\limits_{r\geq 1} \gaussian{n}{r} (-1)^{r-1}q^{\binom{r+1}{2}-rn} \frac{(yq;q)_r}{(zq;q)_r}\sum\limits_{i\geq 1} \gaussian{r}{i} \left(1-z^iq^{(-k+1)i}\frac{(ytz^{-1}q^{k};q)_i}{(ytq;q)_i}\right)\frac{(-1)^{i-1}q^{\binom{i}{2}+ik}}{(1-yq^i)^{k}}
\end{multline}

\begin{multline} \notag
= \frac{1}{(zq^{-k+1};q)_k}\sum\limits_{i \geq 1}  \left(1-z^iq^{(-k+1)i}\frac{(ytz^{-1}q^{k};q)_i}{(ytq;q)_i}\right)\frac{(-1)^{i-1}q^{\binom{i}{2}+ik}}{(1-yq^i)^{k}}\sum\limits_{r\geq 1} \gaussian{n}{r}\gaussian{r}{i} (-1)^{r-1}q^{\binom{r+1}{2}-rn} \frac{(yq;q)_r}{(zq;q)_r}.
\end{multline}
Now using \eqref{lesserbinomtransform3} and making further manipulations, we finally obtain

\begin{multline} \notag
 \sum\limits_{n \geq r_1 \geq r_2 \geq \cdots \geq r_{k-1} \geq r_k \geq 1}\left(1-y^{r_{k}}q^{r_k}\frac{(t;q)_{r_k}}{(ytq;q)_{r_k}}\right)\frac{y^{-r_{k}}(yq;q)_{r_k-1}}{(q;q)_{r_k}} \prod\limits_{j=1}^{k-1} \frac{q^{r_j}}{(1-zq^{r_j-j})(1-yq^{r_j})} 
\end{multline}

\begin{multline} \notag
= \frac{(yq;q)_n}{(q;q)_n(zq;q)_{n-1}(zq^{-k+1};q)_k}\sum\limits_{i \geq 1} \gaussian{n}{i} \left(1-z^iq^{(-k+1)i}\frac{(ytz^{-1}q^{k};q)_i}{(ytq;q)_i}\right)\frac{y^{-i}q^{i(k-1)}(yq;q)_i(zy^{-1};q)_{n-i}}{(1-yq^i)^{k}},
\end{multline}
and the proof is complete.
\end{proof}
Next, let us deploy the Pochhammer symbol $(x)_n = x(x+1)...(x+n-1)$, $(x)_0 =1 $. Now in \textbf{Proposition 3}, multiply both sides by $(1-y)(1-q)^{2k-2}$, put $z= q^{b}, y = q^{a}$, and $t = q^{c}$, and finally let $q \to 1$, to obtain
\vspace{0.7cm}
\begin{corollary} For natural numbers $n$ and $k$, and for complex values $a, b,$ and $c,$ there holds

\begin{multline}
       \tag{4.7} \sum\limits_{n \geq r_1 \geq r_2 \geq \cdots \geq r_{k-1} \geq r_k \geq 1}  \left( \frac{(c)_{r_k}}{(a+c+1)_{r_k}} - \frac{(b-a-k)_{r_k}}{(b-k+1)_{r_k}}\right)\frac{(a)_{r_k}}{(r_k)!} \prod\limits_{j=1}^{k-1} \frac{1}{(b+r_j-j)(a+r_j)} \\ = \frac{(a)_{n+1}}{n! (b-k+1)_{n+k-1}} \sum\limits_{i \geq 1} \binom{n}{i} \frac{(a+1)_i(a+c-b+k)_i (b-a)_{n-i}}{(a+i)^k(a+c+1)_i}.
    \end{multline}

    Provided that $-a \notin \mathbb{N}$, $-a-c \notin \mathbb{N}$, and $b \notin \{k-1, k-2,..., -n+2, -n+1\}$. $\binom{n}{r} = \frac{n(n-1)...(n-r+1)}{r!}$ is the binomial coefficient.
\end{corollary}

Then, let $y=1$ in \textbf{Proposition 3} to arrive at the following corollary.
\vspace{0.7cm}
\begin{corollary} For natural numbers $n$ and $k$, and for complex values $z,t$, there holds

\begin{multline} \tag{4.8} \label{lessergeneralizedGouZhang}
 \sum\limits_{n \geq r_1 \geq r_2 \geq \cdots \geq r_{k-1} \geq r_k \geq 1}\left(\frac{1}{1-zq^{r_k-k}}-\frac{1}{1-tq^{r_k}}\right) \prod\limits_{j=1}^{k-1}\frac{q^{r_j}}{(1-zq^{r_j-j})(1-q^{r_j})}  \\= \frac{1}{(zq^{-k+1};q)_{k+n-1}}\sum\limits_{i \geq 1} \gaussian{n}{i} \frac{z^i(q;q)_{i-1}(tz^{-1}q^{k};q)_i(z;q)_{n-i}}{(1-q^i)^{k-1}(tq;q)_i}.
\end{multline}

Provided that $z \neq q^{r}$, for $-n+1 \leq r \leq k-1$, $r \in \mathbb{Z}$ and $t \neq q^{-s}$, for $1 \leq s \leq n$, $s \in \mathbb{N}$.
\end{corollary}

Which is a generalization of Guo and Zhang's results, \eqref{gouzhang} and Theorem 4.1 in \cite{GuoZhang}. This identity is also proved by the author in \cite{maw2024qseriesidentitiesmultiplesums}. When we let $n \to \infty$, \textbf{Proposition 3} gives the following infinite $q$-series identity.
\vspace{0.7cm}
\begin{corollary} For natural number $k$ and a complex value $|q|<1$, and for complex values $y,z,t$, such that $|\frac{z}{y}| < 1$, there holds

\begin{multline} \tag{4.9} \label{infinitegeneralizedGouZhang}
 \sum\limits_{r_1 \geq r_2 \geq \cdots \geq r_{k-1} \geq r_k \geq 1}\left(\frac{(t;q)_{r_k}}{(ytq;q)_{r_k}}-\frac{(zy^{-1}q^{-k};q)_{r_k}}{(zq^{-k+1};q)_{r_k}}\right)\frac{q^{r_k}(yq;q)_{r_k-1}}{(q;q)_{r_k}} \prod\limits_{j=1}^{k-1}\frac{q^{r_j}}{(1-zq^{r_j-j})(1-yq^{r_j})}\\= \frac{(yq;q)_\infty(zy^{-1};q)_{\infty}}{(q;q)_\infty(zq^{-k+1};q)_{\infty}}\sum\limits_{i = 1}^{\infty} \frac{z^iy^{-i}(yq;q)_i(ytz^{-1}q^{k};q)_i}{(1-yq^i)^{k}(q;q)_i(ytq;q)_i}.
\end{multline}
\end{corollary}

For $k=1$, \textbf{Proposition 3} gives 
\vspace{0.7cm}
\begin{corollary} For natural number $n$ and complex values $y,z,t$, there holds
    \begin{multline} \tag{4.10} \label{lessergeneralizedGouZhang2}
 \sum\limits_{n \geq r \geq 1} q^r\left(\frac{(t;q)_{r}}{(ytq;q)_{r}}-\frac{(zy^{-1}q^{-1};q)_{r}}{(z;q)_{r}}\right)\frac{(yq;q)_{r-1}}{(q;q)_{r}} = \frac{(yq;q)_n}{(q;q)_n(z;q)_{n}}\sum\limits_{i \geq 1} \gaussian{n}{i} \frac{z^iy^{-i}(yq;q)_i(ytz^{-1}q;q)_i(zy^{-1};q)_{n-i}}{(1-yq^i)(ytq;q)_i}.
\end{multline}
Provided that $z \neq q^{-r}$, $r \in \mathbb{N} \cup \{0\}$, $t \neq y^{-1}q^{-m}$, $m \in \mathbb{N}$.
\end{corollary}

The next corollary is the case $y_1=...=y_k=0$ of \textbf{Corollary 9}.
\vspace{0.7cm}
\begin{corollary} For natural numbers $n$ and $k$, and for complex values $z_1,...,z_k$, there holds

\begin{multline} \notag
      \frac{(q;q)_n}{(z_1q;q)_n}\sum\limits_{n \geq r_1 \geq r_2 \geq \cdots \geq r_{k-1} \geq r_k \geq 1} \frac{(1-(x;q)_{r_k})q^{r_{k}}(z_{k}q;q)_{r_{k}-1}}{(q;q)_{r_{k}}} \prod\limits_{j=1}^{k-1} \frac{q^{r_j}(z_jq;q)_{r_j-1}}{(z_{j+1}q;q)_{r_j}} \\ \tag{4.11} \label{reducedanalogguassainmultiple}
= \sum\limits_{r\geq 1} \gaussian{n}{r} \frac{(-1)^{r-1}x^rq^{\binom{r}{2}+rk}}{(1-z_1q^r)(1-z_2q^r)...(1-z_kq^r)} .
\end{multline}

Provided that no $z_j \neq q^{-m}$, $m \in \mathbb{N}$.
\end{corollary}

Multiply both sides of the above equality by $(1-q)^k$, put $z_i = q^{a_i}$, and let $q \to 1$, to get 
\vspace{0.7cm}
\begin{corollary} For natural numbers $n$ and $k$, and for complex values $a_1,...,a_k$, there holds

\begin{multline}
    \tag{4.12} \frac{n!}{(a_1+1)_n}\sum\limits_{n \geq r_1 \geq r_2 \geq \cdots \geq r_{k-1} \geq r_k \geq 1} \frac{(1-(1-x)^{r_k})(a_k+1)_{r_k-1}}{r_k!} \prod\limits_{j=1}^{k-1} \frac{(a_j+1)_{r_j-1}}{(a_{j+1}+1)_{r_j}} \\ = \sum\limits_{r\geq 1} \binom{n}{r} \frac{(-1)^{r-1}x^r}{(a_1+r)(a_2+r)...(a_k+r)}.
\end{multline}

Provided that $-a_j \notin \mathbb{N}$.
\end{corollary}

Now we end with the following case. We use \textbf{Lemma 4} on \eqref{lesserbinomtransform2}, and replace $z$ by $z^{-1}$, and $y$ by  $y^{-1}$ to get

\[ \tag{4.13} \label{lesserbinomtransform2alt}
\sum\limits_{r \geq 1} \gaussian{n}{r} (-1)^{r-1}q^{\binom{r}{2}}x^ry^r \frac{(z;q)_r}{(yz;q)_r} = \frac{1}{(yz;q)_n} \sum\limits_{r \geq 1} \gaussian{n}{r} (1-(x;q)_r)y^r(z;q)_r(y;q)_{n-r}.
\]

Now apply \textbf{Theorem 1} to the above equality with $A_r(q) = \gaussian{n}{r} (-1)^{r-1}q^{\binom{r}{2}}y^r \frac{(z;q)_r}{(yz;q)_r}$ and $B_r(q) = \gaussian{n}{r}\frac{y^r(z;q)_r(y;q)_{n-r}}{(yz;q)_n}$, then we arrive at the following proposition.
\vspace{0.7cm}
\begin{proposition}
For a natural number $n$, for complex values $z,y,x$, and for non-negative integers, $m_1, m_2,...,m_k$ and $n_1,n_2,...,n_k$, there holds
\begin{multline} \tag{4.14} \label{analogmultiplegaussianbinomtransform}
\sum\limits_{r \geq 1} \gaussian{n}{r} (-1)^{r-1}q^{\binom{r}{2}}y^r \frac{(z;q)_r}{(yz;q)_r} H_r[m_1+1, \{1\}^{n_1-1}, m_2+1, \{1\}^{n_2-1}, m_3+1,..., \{1\}^{n_{k-1}-1}, m_k+1, \{1\}^{n_k} ; x] \\= \frac{1}{(yz;q)_n} \sum\limits_{r \geq 1} \gaussian{n}{r}y^r(z;q)_r(y;q)_{n-r} U_r[\{1\}^{m_1}, n_1+1, \{1\}^{m_2-1}, n_2+1, \{1\}^{m_3-1}, n_3+1,..., \{1\}^{m_k-1},n_k+1;x].
\end{multline}
 Provided that $yz \neq q^{-r}, r \in \mathbb{N} \cup \{0\}$.
\end{proposition}
\vspace{0.7cm}
\section{Certain basic hypergeometric multiple sums}
In this section, we proceed to examine certain kinds of basic hypergeometric multiple sums. For our purpose, we will first need to grasp the use of the following operators. Let the operators $K_z$ and $L_{z;t}$ be defined as follows.

\[
K_z = \frac{(z;q)_\infty}{(q;q)_\infty} \sum\limits_{j \geq 0} \frac{(-1)^j q^{\binom{j+1}{2}}}{(q;q)_j}\eta_z^j,
\]

\[
L_{z;t} = \frac{(t;q)_\infty}{(zt;q)_\infty} \sum\limits_{j \geq 0} \frac{(z;q)_j}{(q;q)_j}t^j \eta_z^j.
\]

And from our proof of \textbf{Lemma 12}, the following transformations are evident.

\[
K_z \hspace{0.025cm} \frac{1}{1-zq^n} = (z;q)_n,
\]
\[
K_z \hspace{0.025cm} z^{n-r}\frac{(zq;q)_{r-1}}{(zq;q)_n} = z^{n-r}\frac{(z;q)_r}{(q;q)_{n-r}},
\]

\[
L_{z;t} \hspace{0.025cm} z^n = z^n\frac{(t;q)_n}{(zt;q)_n},
\]

\[
L_{z;t} \hspace{0.025cm} (z;q)_n = \frac{(z;q)_n}{(zt;q)_n},
\]

\[
L_{z;t} \hspace{0.025cm} z^{n-r} (z;q)_r = z^{n-r}\frac{(z;q)_r(t;q)_{n-r}}{(zt;q)_n}.
\]

Then if we define the operator $G_{z;t} = L_{z;t} K_z$, we have the transformations 

\[ \tag{5.1} \label{Gtransformation1}
G_{z;t} \hspace{0.025cm} \frac{1}{1-zq^n} = \frac{(z;q)_n}{(zt;q)_n},
\]
\[ \tag{5.2} \label{Gtransformation2}
G_{z;t} \hspace{0.025cm} z^{n-r}\frac{(zq;q)_{r-1}}{(zq;q)_n} = z^{n-r}\frac{(z;q)_r(t;q)_{n-r}}{(zt;q)_n(q;q)_{n-r}}.
\]

Utilising the operator $G$, we shall prove the following basic hypergeometric multiple sum identity.

\vspace{0.7cm}

\begin{proposition}[A basic hypergeometric multiple sum identity] Let $A_1(q), A_2(q), ...$ and $B_1(q), B_2(q), ...$ be sequences not depending on $x$, satisfying the relation 

\[
\sum\limits_{r\geq 1} A_r(q) x^r = \sum\limits_{r \geq 1} B_r(q) (1-(x;q)_r),
\]

for all complex values $x$, then

\begin{multline} \notag
 \sum\limits_{r\geq 1} A_r(q) \frac{(q;q)_r}{(z_1w_1;q)_r}\sum\limits_{r_0=r \geq r_1 \geq r_2 \geq \cdots \geq r_{k-1} \geq r_k \geq 1}\frac{x^{r_k}t_k^{r_{k}}z_k^{r_{k-1}-r_k}(y_k;q)_{r_k}(z_{k};q)_{r_{k}}(w_k;q)_{r_{k-1}-r_k}}{(y_kt_k;q)_{r_k}(q;q)_{r_{k}}(q;q)_{r_{k-1}-r_{k}}} \times \\ \times \prod\limits_{j=1}^{k-1} \frac{t_j^{r_j}z_j^{r_{j-1}-r_j}(y_j;q)_{r_j}(z_j;q)_{r_j}(w_j;q)_{r_{j-1}-r_j}}{(y_jt_j;q)_{r_j}(z_{j+1}w_{j+1};q)_{r_j}(q;q)_{r_{j-1}-r_j}}
\end{multline} 

\begin{multline} 
    = \sum\limits_{r\geq 1} B_r(q) \frac{(z_1;q)_r(q;q)_r}{(z_1w_1;q)_r(y_1t_1;q)_r}\sum\limits_{r_0=r \geq r_1 \geq r_2 \geq \cdots \geq r_{k-1} \geq r_k \geq 1}\frac{(1-(x;q)_{r_k})t_k^{r_k}(y_k;q)_{r_k}(t_k;q)_{r_{k-1}-r_k}}{(q;q)_{r_k}(q;q)_{r_{k-1}-r_k}} \times \\ \tag{5.3} \label{basichypergeometricmultiplesum} \times \prod\limits_{j=1}^{k-1} \frac{t_j^{r_j}(z_{j+1};q)_{r_j}(y_j;q)_{r_j}(t_j;q)_{r_{j-1}-r_j}}{(z_{j+1}w_{j+1};q)_{r_j}(y_{j+1}t_{j+1};q)_{r_j}(q;q)_{r_{j-1}-r_j}}, 
\end{multline}

for complex values $z_1,z_2,...,z_k$, $y_1,y_2,...,y_k$, $w_1,w_2,...,w_k$, and $t_1,t_2,...,t_k$, provided that $z_jw_j, y_jt_j \neq q^{-r}, r \in \mathbb{N} \cup \{0\}$.
\end{proposition}

Now we present our proof of \textbf{Proposition 21}.
\begin{proof}[Proof of Proposition 21]

Let $A_1(q), A_2(q), ...$ and $B_1(q), B_2(q), ...$ be arbitrary sequences not depending on $x$, satisfying the relation 

\[
\sum\limits_{r\geq 1} A_r(q) x^r = \sum\limits_{r \geq 1} B_r(q) (1-(x;q)_r),
\]

for all arbitrary $x$. Then by \textbf{Lemma 4}, we have

\[
\sum\limits_{r\geq 1} B_r(q^{-1}) x^r = \sum\limits_{r \geq 1} A_r(q^{-1}) (1-(x;q)_r).
\]

Applying \textbf{Theorem 2} to the above equality, we arrive at

\begin{multline} \notag
    \sum\limits_{r\geq 1} B_r(q^{-1}) \frac{q^r(q;q)_r}{(1-z_1q^r)(y_1q;q)_r}\sum\limits_{r_0=r \geq r_1 \geq r_2 \cdots \geq r_{k-1} \geq r_k \geq 1} \frac{x^{r_{k}}y_{k}^{r_{k-1}-r_{k}}(y_kq;q)_{r_k-1}}{(q;q)_{r_k}} \prod\limits_{j=1}^{k-1} \frac{q^{r_j}y_j^{r_{j-1}-r_j}(y_jq;q)_{r_{j}-1}}{(1-z_{j+1}q^{r_j})(y_{j+1}q;q)_{r_j}}
\end{multline}

\begin{multline} \notag
= \sum\limits_{r\geq 1} A_r(q^{-1}) \frac{(q;q)_r}{(z_1q;q)_r}\sum\limits_{r_0=r \geq r_1 \geq r_2 \geq \cdots \geq r_{k-1} \geq r_k \geq 1}\frac{(1-(x;q)_{r_k})q^{r_{k}}(z_{k}q;q)_{r_{k}-1}}{(1-y_{k}q^{r_{k}})(q;q)_{r_{k}}} \prod\limits_{j=1}^{k-1} \frac{q^{r_j}(z_jq;q)_{r_j-1}}{(1-y_jq^{r_j})(z_{j+1}q;q)_{r_j}},
\end{multline}

for all arbitrary values $z_1,z_2,...,z_k$ and $y_1,y_2,...,y_k$. Next, we apply the series of operators $G_{y_k;t_k}...G_{y_2;t_2}G_{y_1;t_1}$ to both sides of the above equality. Then by \eqref{Gtransformation1} and \eqref{Gtransformation2}, we obtain

\begin{multline} \notag
    \sum\limits_{r\geq 1} B_r(q^{-1}) \frac{q^r(q;q)_r}{(1-z_1q^r)(y_1t_1;q)_r}\sum\limits_{r_0=r \geq r_1 \geq r_1 \geq \cdots \geq r_{k-1} \geq r_k \geq 1}\frac{x^{r_{k}}y_{k}^{r_{k-1}-r_{k}}(y_k;q)_{r_k}(t_k;q)_{r_{k-1}-r_k}}{(q;q)_{r_k}(q;q)_{r_{k-1}-r_k}}  \times \\ \times \prod\limits_{j=1}^{k-1}\frac{q^{r_j}y_j^{r_{j-1}-r_j}(y_j;q)_{r_j}(t_j;q)_{r_{j-1}-r_j}}{(1-z_{j+1}q^{r_j})(y_{j+1}t_{j+1};q)_{r_j}(q;q)_{r_{j-1}-r_j}}
\end{multline}

\begin{multline} \notag
= \sum\limits_{r\geq 1} A_r(q^{-1}) \frac{(q;q)_r}{(z_1q;q)_r}\sum\limits_{r_0=r \geq r_1 \geq r_2 \geq \cdots \geq r_{k-1} \geq r_k \geq 1}\frac{(1-(x;q)_{r_k})q^{r_{k}}(y_k;q)_{r_k}(z_{k}q;q)_{r_{k}-1}}{(y_kt_k;q)_{r_k}(q;q)_{r_{k}}} \prod\limits_{j=1}^{k-1}\frac{q^{r_j}(y_j;q)_{r_j}(z_jq;q)_{r_j-1}}{(y_jt_j;q)_{r_j}(z_{j+1}q;q)_{r_j}}.
\end{multline}

Now we again use \textbf{Lemma 4} on the above equality, replacing $z_i$ by $ z_i^{-1}$, $y_i$ by $ y_i^{-1}$, and $t_i$ by $t_i^{-1}$, for $1 \leq i \leq k$, to yield

\begin{multline} \notag
    \sum\limits_{r\geq 1} B_r(q) \frac{(q;q)_r}{(1-z_1q^r)(y_1t_1;q)_r}\sum\limits_{r_0=r \geq r_1 \geq r_2 \geq \cdots \geq r_{k-1} \geq r_k \geq 1} \frac{(1-(x;q)_{r_k})t_k^{r_k}(y_k;q)_{r_k}(t_k;q)_{r_{k-1}-r_k}}{(q;q)_{r_k}(q;q)_{r_{k-1}-r_k}} \times \\ \times \prod\limits_{j=1}^{k-1} \frac{t_j^{r_j}(y_j;q)_{r_j}(t_j;q)_{r_{j-1}-r_j}}{(1-z_{j+1}q^{r_j})(y_{j+1}t_{j+1};q)_{r_j}(q;q)_{r_{j-1}-r_j}}
\end{multline}

\begin{multline} \notag
= \sum\limits_{r\geq 1} A_r(q) \frac{(q;q)_r}{(z_1q;q)_r}\sum\limits_{r_0=r \geq r_1 \geq r_2 \geq \cdots \geq r_{k-1} \geq r_{k} \geq 1} \frac{x^{r_k}t_k^{r_{k}}z_k^{r_{k-1}-r_k}(y_k;q)_{r_k}(z_{k}q;q)_{r_{k}-1}}{(y_kt_k;q)_{r_k}(q;q)_{r_{k}}} \prod\limits_{j=1}^{k-1} \frac{t_j^{r_j}z_j^{r_{j-1}-r_j}(y_j;q)_{r_j}(z_jq;q)_{r_j-1}}{(y_jt_j;q)_{r_j}(z_{j+1}q;q)_{r_j}}.
\end{multline}

Now, we again apply the series of operators $G_{z_k;w_k}...G_{z_2;w_2}G_{z_1;w_1}$ to both sides of the above equality. By \eqref{Gtransformation1} and \eqref{Gtransformation2}, we finally obtain

\begin{multline} \notag
    \sum\limits_{r\geq 1} B_r(q) \frac{(z_1;q)_r(q;q)_r}{(z_1w_1;q)_r(y_1t_1;q)_r}\sum\limits_{r_0=r \geq r_1 \geq r_2 \geq \cdots \geq r_{k-1} \geq r_k \geq 1} \frac{(1-(x;q)_{r_k})t_k^{r_k}(y_k;q)_{r_k}(t_k;q)_{r_{k-1}-r_k}}{(q;q)_{r_k}(q;q)_{r_{k-1}-r_k}} \times \\ \times \prod\limits_{j=1}^{k-1}\frac{t_j^{r_j}(z_{j+1};q)_{r_j}(y_j;q)_{r_j}(t_j;q)_{r_{j-1}-r_j}}{(z_{j+1}w_{j+1};q)_{r_j}(y_{j+1}t_{j+1};q)_{r_j}(q;q)_{r_{j-1}-r_j}}
\end{multline}

\begin{multline} \notag
= \sum\limits_{r\geq 1} A_r(q) \frac{(q;q)_r}{(z_1w_1;q)_r}\sum\limits_{r_0 =r \geq r_1 \geq r_2 \geq \cdots \geq r_{k-1} \geq r_k \geq 1}\frac{x^{r_k}t_k^{r_{k}}z_k^{r_{k-1}-r_k}(y_k;q)_{r_k}(z_{k};q)_{r_{k}}(w_k;q)_{r_{k-1}-r_k}}{(y_kt_k;q)_{r_k}(q;q)_{r_{k}}(q;q)_{r_{k-1}-r_{k}}} \times \\ \times \prod\limits_{j=1}^{k-1} \frac{t_j^{r_j}z_j^{r_{j-1}-r_j}(y_j;q)_{r_j}(z_j;q)_{r_j}(w_j;q)_{r_{j-1}-r_j}}{(y_jt_j;q)_{r_j}(z_{j+1}w_{j+1};q)_{r_j}(q;q)_{r_{j-1}-r_j}}.
\end{multline} 

\end{proof}

We use \textbf{Proposition 21} on $q$-binomial theorem \eqref{gaussbinom}, with $A_r(q) = \gaussian{n}{r} (-1)^{r-1}q^{\binom{r}{2}}$ and with $B_n(q) =1$, and $B_r(q) = 0$ for all $r \neq n$. Then we have the following corollary.
\vspace{0.7cm}
\begin{corollary} For complex values $x$, $z_1,z_2,...,z_k$, $y_1,y_2,...,y_k$, $w_1,w_2,...,w_k$, and $t_1,t_2,...,t_k$, there holds

\begin{multline} \notag
 \sum\limits_{r\geq 1} \gaussian{n}{r} (-1)^{r-1}q^{\binom{r}{2}} \frac{(q;q)_r}{(z_1w_1;q)_r}\sum\limits_{r_0=r \geq r_1 \geq r_2 \geq \cdots \geq r_{k-1} \geq r_k \geq 1}\frac{x^{r_k}t_k^{r_{k}}z_k^{r_{k-1}-r_k}(y_k;q)_{r_k}(z_{k};q)_{r_{k}}(w_k;q)_{r_{k-1}-r_k}}{(y_kt_k;q)_{r_k}(q;q)_{r_{k}}(q;q)_{r_{k-1}-r_{k}}} \times \\ \times \prod\limits_{j=1}^{k-1} \frac{t_j^{r_j}z_j^{r_{j-1}-r_j}(y_j;q)_{r_j}(z_j;q)_{r_j}(w_j;q)_{r_{j-1}-r_j}}{(y_jt_j;q)_{r_j}(z_{j+1}w_{j+1};q)_{r_j}(q;q)_{r_{j-1}-r_j}}
\end{multline} 

\begin{multline} 
    = \frac{(z_1;q)_n(q;q)_n}{(z_1w_1;q)_n(y_1t_1;q)_n}\sum\limits_{r_0=n \geq r_1 \geq r_2 \geq \cdots \geq r_{k-1} \geq r_k \geq 1} \frac{(1-(x;q)_{r_k})t_k^{r_k}(y_k;q)_{r_k}(t_k;q)_{r_{k-1}-r_k}}{(q;q)_{r_k}(q;q)_{r_{k-1}-r_k}} \times \\ \times \prod\limits_{j=1}^{k-1} \frac{t_j^{r_j}(z_{j+1};q)_{r_j}(y_j;q)_{r_j}(t_j;q)_{r_{j-1}-r_j}}{(z_{j+1}w_{j+1};q)_{r_j}(y_{j+1}t_{j+1};q)_{r_j}(q;q)_{r_{j-1}-r_j}}. \tag{5.4} \label{gaussianbasichypergeometricmultiplesum}
\end{multline}

Provided that $z_jw_j, y_jt_j \neq q^{-r}, r \in \mathbb{N} \cup \{0\}$.
\end{corollary}

With $z_1=...=z_k=0$, we arrive at
\vspace{0.7cm}
\begin{corollary}  For complex values $x$, $y_1,y_2,...,y_k$, and $t_1,t_2,...,t_k$, there holds

\begin{multline} 
    \frac{(q;q)_n}{(y_1t_1;q)_n}\sum\limits_{r_0=n \geq r_1 \geq r_2 \geq \cdots \geq r_{k-1} \geq r_k \geq 1} \frac{(1-(x;q)_{r_k})t_k^{r_k}(y_k;q)_{r_k}(t_k;q)_{r_{k-1}-r_k}}{(q;q)_{r_k}(q;q)_{r_{k-1}-r_k}} \prod\limits_{j=1}^{k-1} \frac{t_j^{r_j}(y_j;q)_{r_j}(t_j;q)_{r_{j-1}-r_j}}{(y_{j+1}t_{j+1};q)_{r_j}(q;q)_{r_{j-1}-r_j}} 
   \\ = \sum\limits_{r\geq 1} \gaussian{n}{r} (-1)^{r-1}q^{\binom{r}{2}}(xt_1t_2...t_k)^{r} \frac{(y_1;q)_r(y_2;q)_r...(y_k;q)_r}{(y_1t_1;q)_r(y_2t_2;q)_r...(y_kt_k;q)_r}.
   \tag{5.5} \label{lessergaussianbasichypergeometricmultiplesum}
\end{multline}

Provided that $ y_jt_j \neq q^{-r}, r \in \mathbb{N} \cup \{0\}$.
\end{corollary}

For $n \to \infty$, this gives
\vspace{0.7cm}
\begin{corollary}  For $k \geq 2$ and for complex values $x$, $y_1,y_2,...,y_k$, and $t_1,t_2,...,t_k$, there holds

 \begin{multline} 
    \frac{(t_1;q)_\infty}{(y_1t_1;q)_\infty}\sum\limits_{ r = 1}^{\infty}\frac{t_1^{r}(y_1;q)_{r}}{(y_2t_2;q)_{r}}\sum\limits_{r_0=r \geq r_1 \geq r_2 \geq \cdots \geq r_{k-2} \geq r_{k-1} \geq 1}\frac{(1-(x;q)_{r_{k-1}})t_k^{r_{k-1}}(y_k;q)_{r_{k-1}}(t_k;q)_{r_{k-2}-r_{k-1}}}{(q;q)_{r_{k-1}}(q;q)_{r_{k-2}-r_{k-1}}} \times \\ \times \prod\limits_{j=1}^{k-2} \frac{t_{j+1}^{r_{j+1}}(y_{j+1};q)_{r_{j+1}}(t_{j+1};q)_{r_j-r_{j+1}}}{(y_{j+2}t_{j+2};q)_{r_{j+1}}(q;q)_{r_j-r_{j+1}}}  = \sum\limits_{r= 1}^{\infty}  \frac{(-1)^{r-1}q^{\binom{r}{2}}(xt_1t_2...t_k)^{r}(y_1;q)_r(y_2;q)_r...(y_k;q)_r}{(q;q)_r(y_1t_1;q)_r(y_2t_2;q)_r...(y_kt_k;q)_r}.
   \tag{5.6} \label{infintelessergaussianbasichypergeometricmultiplesum}
\end{multline}

Provided that $|t_j|<1$, and $ y_jt_j \neq q^{-r}, r \in \mathbb{N} \cup \{0\}$.
\end{corollary}

Now if we let $t_1=...=t_k=1$ in \eqref{gaussianbasichypergeometricmultiplesum}, we arrive to the following fact. 
\vspace{0.7cm}
\begin{corollary}  For complex values $x$, $z_1,z_2,...,z_k$, and $w_1,w_2,...,w_k$, there holds

\begin{multline} 
 \sum\limits_{r\geq 1} \gaussian{n}{r} (-1)^{r-1}q^{\binom{r}{2}} \frac{(q;q)_r}{(z_1w_1;q)_r}\sum\limits_{r_0 =r \geq r_1 \geq r_2 \geq \cdots \geq r_{k-1} \geq r_k \geq 1}\frac{(1-x^{r_k})z_k^{r_{k-1}-r_k}(z_{k};q)_{r_{k}}(w_k;q)_{r_{k-1}-r_k}}{(q;q)_{r_{k}}(q;q)_{r_{k-1}-r_{k}}} \times \\ \times \prod\limits_{j=1}^{k-1} \frac{z_j^{r_{j-1}-r_j}(z_j;q)_{r_j}(w_j;q)_{r_{j-1}-r_j}}{(z_{j+1}w_{j+1};q)_{r_j}(q;q)_{r_{j-1}-r_j}}
   = \frac{(x;q)_n(z_1;q)_n(z_2;q)_n...(z_k;q)_n}{(z_1w_1;q)_n(z_2w_2;q)_n...(z_kw_k;q)_n}.
    \tag{5.7} \label{lessergaussianbasichypergeometricmultiplesum2}
\end{multline}

Provided that $z_jw_j \neq q^{-r}, r \in \mathbb{N} \cup \{0\}$.
\end{corollary}

Next, let $n \to \infty$ in \eqref{lessergaussianbasichypergeometricmultiplesum2} to arrive at
\vspace{0.7cm}
\begin{corollary}  For complex values $x$, $z_1,z_2,...,z_k$, and $w_1,w_2,...,w_k$, there holds

 \begin{multline} 
 \sum\limits_{r=1}^{\infty}  \frac{(-1)^{r-1}q^{\binom{r}{2}}}{(z_1w_1;q)_r}\sum\limits_{r_0=r \geq r_1 \geq r_2 \geq \cdots \geq r_{k-1} \geq r_k \geq 1} \frac{(1-x^{r_k})z_k^{r_{k-1}-r_k}(z_{k};q)_{r_{k}}(w_k;q)_{r_{k-1}-r_k}}{(q;q)_{r_{k}}(q;q)_{r_{k-1}-r_{k}}} \times \\ \times \prod\limits_{j=1}^{k-1} \frac{z_j^{r_{j-1}-r_j}(z_j;q)_{r_j}(w_j;q)_{r_{j-1}-r_j}}{(z_{j+1}w_{j+1};q)_{r_j}(q;q)_{r_{j-1}-r_j}}
 = \frac{(x;q)_\infty (z_1;q)_\infty (z_2;q)_\infty ...(z_k;q)_\infty}{(z_1w_1;q)_\infty (z_2w_2;q)_\infty ...(z_kw_k;q)_\infty}.
    \tag{5.8} \label{infinitelessergaussianbasichypergeometricmultiplesum2}
\end{multline}

Provided that $z_jw_j \neq q^{-r}, r \in \mathbb{N} \cup \{0\}$.
\end{corollary}
\vspace{0.7cm}

Now let us further investigate sums of the form 

\[
S_m[a_1,a_2,...,a_n; b_0, b_1,..., b_{n-1}] = \sum\limits_{1 \leq i_m \leq ... \leq i_2 \leq i_1 \leq n} a_{i_1}b_{n-i_1}a_{i_2}b_{i_1-i_2}...a_{i_m}b_{i_{m-1}-i_m}.
\]

Then we have the following recurrence relation for $S.$

\begin{flalign*}
    S_m[a_1,a_2,...,a_n; b_0, b_1,..., b_{n-1}] = a_nb_0 S_{m-1}[a_1,a_2,...,a_n; b_0, b_1,..., b_{n-1}]\\+\sum\limits_{1 \leq r_1 < n} a_{r_1}b_{n-r_1}S_{m-1}[a_1,a_2,...,a_{r_1}; b_0, b_1,..., b_{r_1-1}]. \tag{5.9}
\end{flalign*}

By the above recurrence relation, it will be convenient for us to define $S_0[a_1,a_2,...,a_k; b_0, b_1,..., b_{k-1}]=1$. Now, for two arbitrary sequences $a_1,a_2,...,a_n$ and $b_0,b_1,...,b_{n-1}$, and $1\leq k \leq n$ let us define $F_k(z)$ by the formal series

\[
F_k(z) = 1+\sum\limits_{m=1}^{\infty} z^{m} S_m[a_1,a_2,...,a_k; b_0, b_1,..., b_{k-1}]. \tag{5.10} \label{defF}
\]

Then we have

\begin{flalign*}
 F_n(z) &= 1+\sum\limits_{m=1}^{\infty} z^{m} S_m[a_1,a_2,...,a_n; b_0, b_1,..., b_{n-1}]  \\
 &=1+\sum\limits_{m=1}^{\infty} z^{m} \left\{ a_nb_0 S_{m-1}[a_1,a_2,...,a_n; b_0, b_1,..., b_{n-1}]+\sum\limits_{1 \leq r_1 < n} a_{r_1}b_{n-r_1}S_{m-1}[a_1,a_2,...,a_{r_1}; b_0, b_1,..., b_{r_1-1}] \right\}\\
 &=1+ z a_n b_0  F_n(z) + z \sum\limits_{1 \leq r_1 < n} a_{r_1}b_{n-r_1}  F_{r_1}(z).
\end{flalign*}

Therefore, we have the following recurrence relation for $F$

\[
F_n(z) = \frac{1}{1-z a_n b_0}+\frac{z}{1-z a_n b_0} \sum\limits_{1 \leq r_1 < n} a_{r_1}b_{n-r_1}  F_{r_1}(z). \tag{5.11}
\]

Furthermore, since we have 

\[
F_1(z) = 1+\frac{z a_1 b_0}{1- z a_1 b_0}= \frac{1}{1-z a_1 b_0}. 
\]

We have the following representation for $F$.
\vspace{0.7cm}
\begin{proposition} For $F$ defined as in \eqref{defF}, we have
    \begin{flalign*}
        F_n(z) = \frac{1}{1-z a_n b_0}+\frac{z}{1-z a_n b_0} \sum\limits_{1 \leq r_1 < n} \frac{a_{r_1}b_{n-r_1}}{1-z a_{r_1} b_0} + \frac{z^2}{1-z a_n b_0} \sum\limits_{1 \leq < r_2 < r_1 < n} \frac{a_{r_1}b_{n-r_1}a_{r_2}b_{r_1-r_2}}{(1- z a_{r_1} b_0)(1- z a_{r_2} b_0)} \\ +...+ \frac{z^{n-1}}{1-z a_n b_0} \sum\limits_{1 \leq r_{n-1} <...< r_2 < r_1 < n} \frac{a_{r_1}b_{n-r_1}a_{r_2}b_{r_1-r_2}...a_{r_{n-1}}b_{r_{n-2}-r_{n-1}}}{(1- z a_{r_1} b_0)(1- z a_{r_2} b_0)...(1-za_{r_{n-1}}b_0)}. \tag{5.12} \label{doublesymmetric}
    \end{flalign*}
\end{proposition}

Now, application of \textbf{Proposition 27} on \textbf{Corollary 23} gives 
\vspace{0.7cm}
\begin{proposition} Given that the expressions from both sides do not exhibit singularities, there holds
    \begin{multline}
   \frac{(q;q)_n}{(yt;q)_n-wt^n(y;q)_n}\left\{ 1+ w \sum\limits_{1 \leq r_1 < n} \frac{t^{r_1}(y;q)_{r_1}(t;q)_{n-r_1}}{((yt;q)_{r_1}-wt^{r_1}(y;q)_{r_1})(q;q)_{n-r_1}} \right. \\ \left. + w^2 \sum\limits_{1 \leq r_2 < r_1 < n} \frac{t^{r_1+r_2}(y;q)_{r_1}(t;q)_{n-r_1}(y;q)_{r_2}(t;q)_{r_1-r_2}}{((yt;q)_{r_1}-wt^{r_1}(y;q)_{r_1})(q;q)_{n-r_1}((yt;q)_{r_2}-wt^{r_2}(y;q)_{r_2})(q;q)_{r_1-r_2}} + ... \right. \\ \left. ...+ w^{n-1} \sum\limits_{1 \leq r_{n-1} <...<r_2<r_1<n} \ \prod\limits_{\substack{j=1 \\ r_0 = n}}^{n-1}\frac{t^{r_j}(y;q)_{r_j}(t;q)_{r_{j-1}-r_j}}{((yt;q)_{r_j}-wt^{r_j}(y;q)_{r_j})(q;q)_{r_{j-1}-r_j}}  \right\}  \\= \sum\limits_{r\geq 1} \gaussian{n}{r} \frac{(-1)^{r-1}q^{\binom{r}{2}}(1-q^r)(yt;q)_r}{(1-ytq^{r-1})((yt;q)_r-wt^r(y;q)_r)}, \tag{5.13} \label{newexpression1}
\end{multline}
for all natural numbers $n$.
\end{proposition}

\begin{proof}[Proof of Proposition 28] Put $t_1=t_2=...=t_k=t$ and $y_1=y_2=...=y_k=y$ in \textbf{Corollary 23} to get

\begin{multline} 
    \frac{(q;q)_n}{(yt;q)_n}\sum\limits_{r_0=n \geq r_1 \geq r_2 \geq \cdots \geq r_{k-1} \geq r_k \geq 1}\frac{(1-(x;q)_{r_k})t^{r_k}(y;q)_{r_k}(t;q)_{r_{k-1}-r_k}}{(q;q)_{r_k}(q;q)_{r_{k-1}-r_k}} \prod\limits_{j=1}^{k-1}\frac{t^{r_j}(y;q)_{r_j}(t;q)_{r_{j-1}-r_j}}{(yt;q)_{r_j}(q;q)_{r_{j-1}-r_j}}
   \\ = \sum\limits_{r\geq 1} \gaussian{n}{r} (-1)^{r-1}q^{\binom{r}{2}}x^{r}t^{kr} \frac{(y;q)_r^k}{(yt;q)_r^k}.
   \tag{5.14} \label{refinedlessergaussianbasichypergeometricmultiplesum1}
\end{multline}

Some elementary adjustments give 

\begin{multline} 
    \frac{(q;q)_n}{(yt;q)_n}\sum\limits_{r_0=n \geq r_1 \geq r_2 \geq \cdots \geq r_{k-1} \geq r_k \geq 1}\frac{t^{r_k}(y;q)_{r_k}(t;q)_{r_{k-1}-r_k}(x;q)_{r_k}}{(q;q)_{r_k}(q;q)_{r_{k-1}-r_k}} \prod\limits_{j=1}^{k-1}\frac{t^{r_j}(y;q)_{r_j}(t;q)_{r_{j-1}-r_j}}{(yt;q)_{r_j}(q;q)_{r_{j-1}-r_j}}
   \\ = \sum\limits_{r\geq 1} \gaussian{n}{r} (-1)^{r-1}q^{\binom{r}{2}}(1-x^{r})t^{kr} \frac{(y;q)_r^k}{(yt;q)_r^k}.
   \tag{5.15} \label{refinedlessergaussianbasichypergeometricmultiplesum2}
\end{multline}

Now we apply \textbf{Lemma 5} with the introduction of a new parameter $z$.

\begin{multline} 
    \frac{(q;q)_n}{(yt;q)_n}\sum\limits_{r_0=n \geq r_1 \geq r_2 \geq \cdots \geq r_{k-1} \geq r_k \geq 1}\frac{t^{r_k}(y;q)_{r_k}(t;q)_{r_{k-1}-r_k}(x;q)_{r_k}}{(q;q)_{r_k}(q;q)_{r_{k-1}-r_k}(xz;q)_{r_k}} \prod\limits_{j=1}^{k-1}\frac{t^{r_j}(y;q)_{r_j}(t;q)_{r_{j-1}-r_j}}{(yt;q)_{r_j}(q;q)_{r_{j-1}-r_j}}
   \\ = \sum\limits_{r\geq 1} \gaussian{n}{r} (-1)^{r-1}q^{\binom{r}{2}}(1-x^{r}\frac{(z;q)_r}{(xz;q)_r})t^{kr} \frac{(y;q)_r^k}{(yt;q)_r^k}.
   \tag{5.16} \label{refinedlessergaussianbasichypergeometricmultiplesum3}
\end{multline}

Put $x=q, z=ytq^{-1}$ to arrive at

\begin{multline} 
    \frac{(q;q)_n}{(yt;q)_n}\sum\limits_{r_0=n \geq r_1 \geq r_2 \geq \cdots \geq r_{k-1} \geq r_k \geq 1} \prod\limits_{j=1}^{k}\frac{t^{r_j}(y;q)_{r_j}(t;q)_{r_{j-1}-r_j}}{(yt;q)_{r_j}(q;q)_{r_{j-1}-r_j}} = \sum\limits_{r\geq 1} \gaussian{n}{r} (-1)^{r-1}q^{\binom{r}{2}}\left(\frac{1-q^r}{1-ytq^{r-1}} \right)t^{kr} \frac{(y;q)_r^k}{(yt;q)_r^k}.
   \tag{5.17} \label{refinedlessergaussianbasichypergeometricmultiplesum3b}
\end{multline}

Application of \textbf{Proposition 27} with the introduction of a new parameter $w$, now gives \eqref{newexpression1}.

\end{proof}

Similarly, application of \textbf{Proposition 27} on \textbf{Corollary 25} gives 
\vspace{0.7cm}

\begin{proposition} Given that the expressions from both sides do not exhibit singularities, there holds
    \begin{multline}
   \sum\limits_{r\geq 1} \gaussian{n}{r} (-1)^{r-1}q^{\binom{r+1}{2}-nr}\frac{(q;q)_r}{(zw;q)_r-tw^r(z;q)_r} \left\{1+ t \sum\limits_{1 \leq i_1 < r} \frac{w^{i_1}(z;q)_{i_1}(w;q)_{r-i_1}}{((zw;q)_{i_1}-tw^{i_1} (z;q)_{i_1})(q;q)_{r-i_1} } \right. \\ \left. + t^2\sum\limits_{1 \leq i_2 <i_1 < r} \frac{w^{i_1+i_2}(z;q)_{i_1}(w;q)_{r-i_1}(z;q)_{i_2}(w;q)_{i_1-i_2}}{((zw;q)_{i_1}-tw^{i_1} (z;q)_{i_1})(q;q)_{r-i_1}((zw;q)_{i_2}-tw^{i_2} (z;q)_{i_2})(q;q)_{i_1-i_2} } + ... \right. \\ \left. ... + t^{r-1}\sum\limits_{1 \leq  i_{r-1}< ... <i_2 <i_1 < r} \  \prod\limits_{\substack{j=1 \\ i_0 = r}}^{r-1}  \frac{w^{i_j}(z;q)_{i_j}(w;q)_{i_{j-1}-i_j}}{((zw;q)_{i_j}-tw^{i_j} (z;q)_{i_j})(q;q)_{i_{j-1}-i_j} } \right\}\\ = \left(\frac{1-q^n}{1-zwq^{n-1}}\right)\frac{(zw;q)_n}{(zw;q)_n-tw^n (z;q)_n} , \tag{5.18} \label{newexpression2}
\end{multline}
for all natural numbers $n$.
\end{proposition}

\begin{proof}[Proof of Proposition 29]

We let $z_1=z_2=...=z_{k}=z, w_1=w_2=...=w_k=w$ in \textbf{Corollary 25} to get

\begin{multline} 
 \sum\limits_{r\geq 1} \gaussian{n}{r} (-1)^{r-1}q^{\binom{r}{2}}z^r \frac{(q;q)_r}{(zw;q)_r}\sum\limits_{r_0=r \geq r_1 \geq r_2 \geq \cdots \geq r_{k-1} \geq r_k \geq 1}\frac{(1-x^{r_k})z^{-r_k}(z;q)_{r_{k}}(w;q)_{r_{k-1}-r_k}}{(q;q)_{r_{k}}(q;q)_{r_{k-1}-r_{k}}} \times \\ \times \prod\limits_{j=1}^{k-1}\frac{(z;q)_{r_j}(w;q)_{r_{j-1}-r_j}}{(zw;q)_{r_j}(q;q)_{r_{j-1}-r_j}} = \frac{(x;q)_n(z;q)_n^k}{(zw;q)_n^k}.
    \tag{5.19} \label{refinedlessergaussianbasichypergeometricmultiplesum2b}
\end{multline}

Now we use \textbf{Lemma 4} and replace $w$ by $w^{-1}$, $z$ by $z^{-1}$

\begin{multline} 
 \sum\limits_{r\geq 1} \gaussian{n}{r} (-1)^{r-1}q^{\binom{r+1}{2}-nr} \frac{(q;q)_r}{(zw;q)_r}\sum\limits_{r_0=r \geq r_1 \geq r_2 \geq \cdots \geq r_{k-1} \geq r_k \geq 1}\frac{w^{r_k}(x;q)_{r_k}(z;q)_{r_{k}}(w;q)_{r_{k-1}-r_k}}{(q;q)_{r_{k}}(q;q)_{r_{k-1}-r_{k}}} \times \\ \times \prod\limits_{j=1}^{k-1}\frac{w^{r_j}(z;q)_{r_j}(w;q)_{r_{j-1}-r_j}}{(zw;q)_{r_j}(q;q)_{r_{j-1}-r_j}} = \frac{(1-x^n)w^{nk}(z;q)_n^k}{(zw;q)_n^k}.
    \tag{5.20} \label{refined2lessergaussianbasichypergeometricmultiplesum2}
\end{multline}

Next we use \textbf{Lemma 5} with the introduction of a new parameter $y$ to arrive at

\begin{multline} 
\sum\limits_{r\geq 1} \gaussian{n}{r} (-1)^{r-1}q^{\binom{r+1}{2}-nr} \frac{(q;q)_r}{(zw;q)_r}\sum\limits_{r_0=r \geq r_1 \geq r_2 \geq \cdots \geq r_{k-1} \geq r_k \geq 1}\frac{w^{r_k}(x;q)_{r_k}(z;q)_{r_{k}}(w;q)_{r_{k-1}-r_k}}{(xy;q)_{r_k}(q;q)_{r_{k}}(q;q)_{r_{k-1}-r_{k}}} \times \\ \times \prod\limits_{j=1}^{k-1}\frac{w^{r_j}(z;q)_{r_j}(w;q)_{r_{j-1}-r_j}}{(zw;q)_{r_j}(q;q)_{r_{j-1}-r_j}} = \left(1-x^n\frac{(y;q)_n}{(xy;q)_n}\right)\frac{w^{nk}(z;q)_n^k}{(zw;q)_n^k}.
    \tag{5.21} \label{refined3lessergaussianbasichypergeometricmultiplesum2}
\end{multline}

Put $x=q, y=zwq^{-1}$ to get

\begin{multline} 
 \sum\limits_{r\geq 1} \gaussian{n}{r} (-1)^{r-1}q^{\binom{r+1}{2}-nr} \frac{(q;q)_r}{(zw;q)_r}\sum\limits_{r_0=r \geq r_1 \geq r_2 \geq \cdots \geq r_{k-1} \geq r_k \geq 1} \prod\limits_{j=1}^{k} \frac{w^{r_j}(z;q)_{r_j}(w;q)_{r_{j-1}-r_j}}{(zw;q)_{r_j}(q;q)_{r_{j-1}-r_j}}
   \\ = \left(\frac{1-q^n}{1-zwq^{n-1}}\right)\frac{w^{nk}(z;q)_n^k}{(zw;q)_n^k}.
    \tag{5.22} \label{refined4lessergaussianbasichypergeometricmultiplesum2}
\end{multline}

Application of \textbf{Proposition 27} with the introduction of a new parameter $t$, now gives \eqref{newexpression2}.
    
\end{proof}

Now let us define $\mathcal{H}_r(y;q) = \sum\limits_{i=1}^r \frac{1}{1-yq^{i-1}}$. Then in \eqref{newexpression1}, multiply both sides by $1-t$, let $w=1$ and let $t \to 1$ to get

\vspace{0.7cm}

\begin{corollary} Given that the expressions from both sides do not exhibit singularities, there holds
    \begin{multline}
   \frac{(q;q)_n}{(y;q)_n \mathcal{H}_n(y;q)}\left\{ 1+  \sum\limits_{1 \leq r_1 < n} \frac{1}{\mathcal{H}_{r_1}(y;q)(1-q^{n-r_1})} +  \sum\limits_{1 \leq r_2 < r_1 < n} \frac{1}{\mathcal{H}_{r_1}(y;q)\mathcal{H}_{r_2}(y;q)(1-q^{n-r_1})(1-q^{r_1-r_2})} + ... \right. \\ \left. ...+ \sum\limits_{1 \leq r_{n-1} <...<r_2<r_1<n} \ \prod\limits_{\substack{j=1 \\ r_0 = n}}^{n-1}\frac{1}{\mathcal{H}_{r_j}(y;q)(1-q^{r_{j-1}-r_j})}  \right\}  = \sum\limits_{r\geq 1} \gaussian{n}{r} \frac{(-1)^{r-1}q^{\binom{r}{2}}(1-q^r)}{(1-yq^{r-1})\mathcal{H}_r(y;q)}, \tag{5.23} \label{reciprocalq-harmonic1}
\end{multline}
for all natural numbers $n$.
\end{corollary}

Similarly, in \eqref{newexpression2}, multiply both sides by $1-w$, let $t=1$ and let $w \to 1$ to arrive at

\vspace{0.7cm}

\begin{corollary}
Given that the expressions from both sides do not exhibit singularities, there holds
\begin{multline}
   \sum\limits_{r\geq 1} \gaussian{n}{r} (-1)^{r-1}q^{\binom{r+1}{2}-nr}\frac{(q;q)_r}{(z;q)_r \mathcal{H}_r(z;q)} \left\{1+  \sum\limits_{1 \leq i_1 < r} \frac{1}{\mathcal{H}_{i_1}(z;q)(1-q^{r-i_1})} \right. \\ \left. +  \sum\limits_{1 \leq i_2 <i_1 < r} \frac{1}{\mathcal{H}_{i_1}(z;q)\mathcal{H}_{i_2}(z;q)(1-q^{r-i_1})(1-q^{i_1-i_2})} + ... + \sum\limits_{1 \leq  i_{r-1}< ... <i_2 <i_1 < r} \  \prod\limits_{\substack{j=1 \\ i_0 = r}}^{r-1}  \frac{1}{\mathcal{H}_{i_j}(z;q)(1-q^{i_{j-1}-i_j}) } \right\}\\ = \frac{1-q^n}{(1-zq^{n-1})\mathcal{H}_n(z;q)} , \tag{5.24} \label{reciprocalq-harmonic2}
\end{multline}
for all natural numbers $n$.
\end{corollary}

Let us define the generalized harmonic number $\mathcal{H}_n(a) = \sum\limits_{i=1}^{n} \frac{1}{a+i-1}$. Then, divide both sides of \eqref{reciprocalq-harmonic1} by $1-q$ and let $y=q^a$, then put $q \to 1$ to get

\begin{corollary}Given that the expressions from both sides do not exhibit singularities, there holds
    \begin{multline}
   \frac{n!}{(a)_n \mathcal{H}_n(a)}\left\{ 1+  \sum\limits_{1 \leq r_1 < n} \frac{1}{\mathcal{H}_{r_1}(a)(n-r_1)} +  \sum\limits_{1 \leq r_2 < r_1 < n} \frac{1}{\mathcal{H}_{r_1}(a)\mathcal{H}_{r_2}(a)(n-r_1)(r_1-r_2)} + ... \right. \\ \left. ...+ \sum\limits_{1 \leq r_{n-1} <...<r_2<r_1<n} \ \prod\limits_{\substack{j=1 \\ r_0 = n}}^{n-1}\frac{1}{\mathcal{H}_{r_j}(a)(r_{j-1}-r_j)}  \right\}  = \sum\limits_{r\geq 1} \binom{n}{r} \frac{(-1)^{r-1}r}{(a+r-1)\mathcal{H}_r(a)}, \tag{5.25} \label{reciprocalgeneralizedharmonic1}
\end{multline}
for all natural numbers $n$.
\end{corollary}

Similarly, divide both sides of \eqref{reciprocalq-harmonic2} by $1-q$, let $z = q^a$ and put $q \to 1$ to get 

\vspace{0.7cm}

\begin{corollary}
Given that the expressions from both sides do not exhibit singularities, there holds
\begin{multline}
   \sum\limits_{r\geq 1} \binom{n}{r} (-1)^{r-1}\frac{r!}{(a)_r \mathcal{H}_r(a)} \left\{1+  \sum\limits_{1 \leq i_1 < r} \frac{1}{\mathcal{H}_{i_1}(a)(r-i_1)} +  \sum\limits_{1 \leq i_2 <i_1 < r} \frac{1}{\mathcal{H}_{i_1}(a)\mathcal{H}_{i_2}(a)(r-i_1)(i_1-i_2)} + ... \right. \\ \left. ... + \sum\limits_{1 \leq  i_{r-1}< ... <i_2 <i_1 < r} \  \prod\limits_{\substack{j=1 \\ i_0 = r}}^{r-1}  \frac{1}{\mathcal{H}_{i_j}(a)(i_{j-1}-i_j) } \right\} = \frac{n}{(a+n-1)\mathcal{H}_n(a)} , \tag{5.26} \label{reciprocalgeneralizedharmonic2}
\end{multline}
for all natural numbers $n$.
\end{corollary}

Finally, let $a=1$ in \eqref{reciprocalgeneralizedharmonic1} to get the following identity for reciprocal harmonic number. 

\vspace{0.7cm}

\begin{corollary} For all natural numbers $n$, there holds
    \begin{multline}
   \frac{1}{H_n}\left\{ 1+  \sum\limits_{1 \leq r_1 < n} \frac{1}{H_{r_1}(n-r_1)} +  \sum\limits_{1 \leq r_2 < r_1 < n} \frac{1}{H_{r_1}H_{r_2}(n-r_1)(r_1-r_2)} + ... \right. \\ \left. ...+ \sum\limits_{1 \leq r_{n-1} <...<r_2<r_1<n} \ \prod\limits_{\substack{j=1 \\ r_0 = n}}^{n-1}\frac{1}{H_{r_j}(r_{j-1}-r_j)}  \right\}  = \sum\limits_{r\geq 1} \binom{n}{r} \frac{(-1)^{r-1}}{H_r}, \tag{5.27} \label{reciprocalharmonic}
\end{multline}
where $H_n = \sum\limits_{i=1}^n \frac{1}{i}$ is the harmonic number.
\end{corollary}

\vspace{0.7cm}
\textbf{\Large{Disclosures}}\\
\vspace{0.09cm}
The author reports that there are no financial or non-financial competing interests to declare.
\printbibliography

\end{document}